\documentclass[preprint,12pt]{elsarticle}

\usepackage{lineno}
\usepackage{amssymb,amsxtra,amsmath}
\usepackage{latexsym}
\usepackage{ifthen}
\usepackage{graphicx}
\usepackage{calc,fancyhdr}
\usepackage{indentfirst}
\usepackage{parskip,algorithm,algorithmicx}
\usepackage{algpseudocode}
\usepackage{color}
\usepackage{bigints}
\usepackage{subfig}
\usepackage{Article}
\usepackage{color}
\newcommand{\changeb}[1]{\color{blue}#1\normalcolor}
\newcommand{\changer}[1]{\color{red}#1\normalcolor}

\journal{JCP}

\begin{document}

\begin{frontmatter}

\title{The Overlapped Radial Basis Function-Finite Difference (RBF-FD) Method: A Generalization of RBF-FD}

\author{Varun Shankar}
\address{Department of Mathematics, University of Utah, UT, USA}
\ead{vshankar@math.utah.edu}

\begin{abstract}
We present a generalization of the RBF-FD method that computes RBF-FD weights in finite-sized neighborhoods around the centers of RBF-FD stencils by introducing an overlap parameter $\delta \in [0,1]$ such that $\delta=1$ recovers the standard RBF-FD method and $\delta=0$ results in a full decoupling of stencils. We provide experimental evidence to support this generalization, and develop an automatic stabilization procedure based on local Lebesgue functions for the stable selection of stencil weights over a wide range of $\delta$ values. We provide an a priori estimate for the speedup of our method over RBF-FD that serves as a good predictor for the true speedup. We apply our method to parabolic partial differential equations with time-dependent inhomogeneous boundary conditions-- Neumann in 2D, and Dirichlet in 3D. Our results show that our method can achieve as high as a 60x speedup in 3D over existing RBF-FD methods in the task of forming differentiation matrices.
\end{abstract}

\end{frontmatter}

\section{Introduction}
\label{sec:intro}

Radial basis functions (RBFs) are popular building blocks in the development of numerical methods for partial differential equations (PDEs). RBF interpolants have been used as replacements for polynomial interpolants in generating pseudospectral and finite-difference methods~\cite{ShuDing2003,Bayona2010,Davydov2011,Tolstykh2003,CecilQian2004,Wright200699,Chandhini2007,SPLM}. \changeb{Unlike polynomial-based collocation methods, RBF-based methods can naturally handle irregular collocation node layouts and therefore allow for high-order methods on irregular domains}. RBF-based methods are also increasingly used for the solution of PDEs on node sets that are not unisolvent for polynomials, like the sphere $\mathbb{S}^2$~\cite{FlyerWright:2007,FlyerWright:2009,FoL11,FlyerLehto2012} and other general surfaces~\cite{Piret2012, Piret2016,FuselierWright2013,SWFKIJNMF2014,SWFKJSC2014}.

Unfortunately, interpolation matrices formed from the standard RBF basis have historically been beset by ill-conditioning~\cite{Wendland:2004,Fasshauer:2007}. Error curves for both RBF interpolation and RBF methods for PDEs typically level off as the number of nodes is increased. Fortunately, this phenomenon of ``stagnation'' or ``saturation'' of the errors is not a feature of the \emph{approximation space} spanned by RBFs. This fact has been used to develop many stable algorithms for RBFs with a shape parameter, most applicable for Gaussian RBFs~\cite{FoWr,FornbergPiret:2007,FaMC12,FLF,FoLePo13}, at a cost ten to one hundred times (10--100x) that of standard techniques. However, Flyer et al.~\cite{FlyerPHS} recently demonstrated that adding polynomials of degree up to half the dimension of the RBF space \changeb{overcomes the issue of stagnation errors without significantly increasing the overall cost}. This technique, henceforth referred to as \emph{augmented RBF-FD}, has been used to solve hyperbolic and elliptic PDEs in 2D and 3D domains~\cite{FlyerNS,BarnettPHS}. This new technique has made the use of RBFs without shape parameters feasible for the tasking of generating RBF-generated Finite Difference (RBF-FD) formulas.

The new technique introduces its own difficulties. Since convergence rates are dictated purely by the appended polynomial, achieving high orders of accuracy now requires far more points in both 2D and 3D than if one were to simply use RBFs. \changer{For instance, on a quasi-uniform node set, a 5-point stencil in 2D would roughly correspond to a second order method with standard RBF-FD~\cite{Wright200699}; however, when polynomials are appended, a 12-point stencil would be required for the same order. The situation is exacerbated in 3D, with a 30-point stencil being required for a second-order method when polynomials are appended}. Even parallel implementations cannot fully ameliorate the increased cost to attain a specific order of convergence. The increase in stencil sizes necessitates that fewer number of weights be computed in parallel if one wishes to maximize thread memory occupancy.

We present a generalization of the augmented RBF-FD method designed to decrease the total number of stencils for a given node set and stencil size. Our new method is based on observations of error patterns in interpolation with RBFs and RBFs augmented with polynomials. The new \emph{overlapped} RBF-FD method allows for a large reduction in cost in computing differentiation matrices. Further, the overlapped RBF-FD method allows the use of higher-order methods for only a slightly higher total cost than lower-order ones for the same number of degrees of freedom. We also present a novel stabilization procedure based on \emph{local Lebesgue functions} to stabilize our method in certain rare scenarios. The error estimates for RBF-FD (overlapped or otherwise) are then presented in terms of these local Lebesgue functions. In addition, we develop a local sufficiency condition based on the local Lebesgue functions to ensure that eigenvalues of the discrete RBF-FD Laplacian have only negative real parts. Further, we present a complexity analysis of stable algorithms, augmented RBF-FD and overlapped RBF-FD that establishes the asymptotic speedup of our method over augmented RBF-FD, and augmented RBF-FD over stable algorithms.

The remainder of the paper is organized as follows. In the next section, we briefly review augmented RBF interpolation. Section 3 contains a numerical exploration of error distributions in global approximations with RBF interpolation, augmented RBF interpolation, and polynomial least-squares. The observations from this section are used in Section 4 to motivate the development of the overlapped RBF-FD method. We also present in Section 4 a stabilization technique, a summary of existing error estimates, and a complexity analysis of our method. We describe the implicit time-stepping of the spatially-discretized forced heat equation in the presence of boundary conditions in section 5. Section 6 discusses the eigenvalue stability of the overlapped method when approximating the Laplacian and enforcing boundary conditions. In section 7, we use our method to solve the forced heat equation with time-dependent, inhomogeneous boundary conditions in 2D and 3D, and discuss the effect of overlapping on both convergence and speedup. We conclude with a summary of our results and a discussion of future work in Section 8.
\section{Augmented local RBF interpolation}
\label{sec:rbf_review}

Let $\Omega \subseteq \mathbb{R}^d$, and $\phi :\Omega \times\Omega \to \mathbb{R}$ be a (radial) kernel with the property $\phi(\vx,\vy) := \phi(\|\vx-\vy\|)$ for $\vx,\vy\in\Omega$, where $\|\cdot\|$ is the standard Euclidean norm in $\mathbb{R}^d$. Given a set of nodes $X = \{\vx_k\}_{k = 1}^N \subset \Omega$ and a target function $f:\Omega \to \mathbb{R}$ sampled at the nodes in $X$, we select subsets (henceforth referred to as \emph{stencils}) $\{P_k\}_{k=1}^N$ of the set of nodes $X$, where each stencil $P_k$ consists of the $k$\textsuperscript{th} node $\vx_k$ and its $n-1$ nearest neighbors, where $n << N$. Further, with each stencil $P_k$, we also track an index set $\calI_k = \{\calI^k_1,\hdots,\calI^k_n\}$ that contains the \emph{global} indices of its nodes in the set $X$. The nearest neighbors are typically determined in a preprocessing step using a data structure such as a kd-tree. In~\cite{FlyerPHS}, Flyer et al. demonstrate that it is useful to augment the RBF interpolant with polynomials, with the polynomial of a degree $s$ so that it has $M \lesssim \frac{n}{2}$ basis functions. We adopt this approach, and form \emph{local} RBF interpolants on each stencil $P_k$ such that
\begin{align}
s^k_f(\vx)= \sum\limits_{j={\calI^k_1}}^{\calI^k_n} c_j \phi(\|\vx - \vx_j\|) + \sum\limits_{i=1}^M d_i \psi^k_i(\vx),
\label{eq:rbf_poly}
\end{align}
where superscripts index the stencil $P_k$. This can be written as the block linear system:
\begin{align}
\underbrace{
\begin{bmatrix}
A_k & \Psi_k \\
\Psi_k^T & 0
\end{bmatrix}}_{\hat{A}_k}
\underbrace{
\begin{bmatrix}
{\bf c} \\
{\bf d}
\end{bmatrix}}_{\hat{c}_k}
=
\begin{bmatrix}
\vf \\
{\bf 0}
\end{bmatrix},
\end{align}
where $A_k$ is the RBF interpolation matrix on $P_k$ and $\Psi_k = \psi^k_i(\vx_j)$. \changer{If $s$ is the degree of the appended polynomial, we find $s$ so that:
\begin{align}
\binom{s+d}{d} = \frac{n}{2},
\end{align}
where $n$ is the stencil size, $M$ the number of polynomial terms, and $d$ the number of spatial dimensions. If $d=2$, this is a quadratic equation for $s$, from which we select the positive solution $s = \lf \lfloor\frac{1}{2} \lf(\sqrt{4n + 1} - 3 \rt) \rt\rfloor$. For $n=36$, this gives $s = 4$, and $M = 0.5(s+1)(s+2) = 15$. Similarly, for $n=100$, we obtain $s = 8$ and $M=45$. In 3D, we analytically solve a cubic equation for $s$ and select $M$ similarly. This approach allows us to choose a stencil size, and always append a polynomial of appropriate degree. }

These local interpolants can be used to approximate functions to high-order algebraic accuracy determined by the degree of the appended polynomial. When differentiated, this local interpolation approach can be used to generate scattered-node finite difference (FD) formulas, known in the literature as RBF-FD; this will be explained in greater detail in Section 4. It has been shown that $\hat{A}_k$ is invertible if the node set $X$ is unisolvent for the polynomials $\psi^k$~\cite{Fasshauer:2007}. More interestingly, $\hat{A}_k$ appears to be invertible even if the node set is not known to be polynomial unisolvent~\cite{FlyerPHS}. This latter property is not yet fully understood.

In~\cite{FlyerPHS,BarnettPHS}, Flyer et al. compare augmented local RBF interpolation to polynomial least squares, and conclude that the former is more accurate for a given $n$. Further, they demonstrate that using augmented RBFs can help alleviate stagnation errors~\cite{FlyerNS}. They also note that convergence rates now depend purely on the degree of the appended polynomial, a fact rigorously proved recently by Davydov and Schaback~\cite{DavydovMinimal2016,DavydovSchaback2016}, whose results are discussed further in Section 4.

Before we present the overlapped RBF-FD method, we attempt to provide some further insights into the \emph{error distribution} in augmented RBF interpolation. In the next section, we compare augmented RBF interpolation to polynomial least squares by exploring the error distribution when interpolating a 2D analog of the Runge function. The results from this section are used to motivate the new overlapped RBF-FD method.
\section{Error distributions in global approximation schemes}
\label{sec:runge}

In this section, we explore the error distributions for global approximation with RBFs, augmented RBFs, and polynomial least squares on a simple 2D test case. The insights from global interpolation will be applied to local interpolation as well (see Sections 4 and 6.1). We do so by interpolating the 2D Runge function given by $f(x,y) = \frac{1}{1 + \gamma(x^2 + y^2)}$, with $(x,y) \in  [-1,1]^2$; we set $\gamma = 1$, but our results and observations carry to higher values of $\gamma$ as well, albeit with higher errors. While such studies have been previously performed for infinitely-smooth RBFs~\cite{FZ07}, our focus is on polyharmonic splines augmented with polynomials. We use two types of node sets for collocation: Cartesian and Halton nodes. However, we always evaluate all approximants at $10^5$ Cartesian nodes. This large number of evaluation nodes allows us to resolve fine features in the error distributions. We use the heptic RBF given by $\phi(r) = r^7$ for all tests. However, the results here carry over to the other polyharmonic spline RBFs with larger exponents of $r$ leading to lower errors. We denote the number of nodes by $n$ rather than $N$ to emphasize that these tests will be applicable to individual RBF-FD stencils as well.

\subsection{Cartesian nodes}
\label{sec:runge_cart}
\begin{figure}[tbhp]
\centering
\subfloat[]
{
	\includegraphics[scale=0.4]{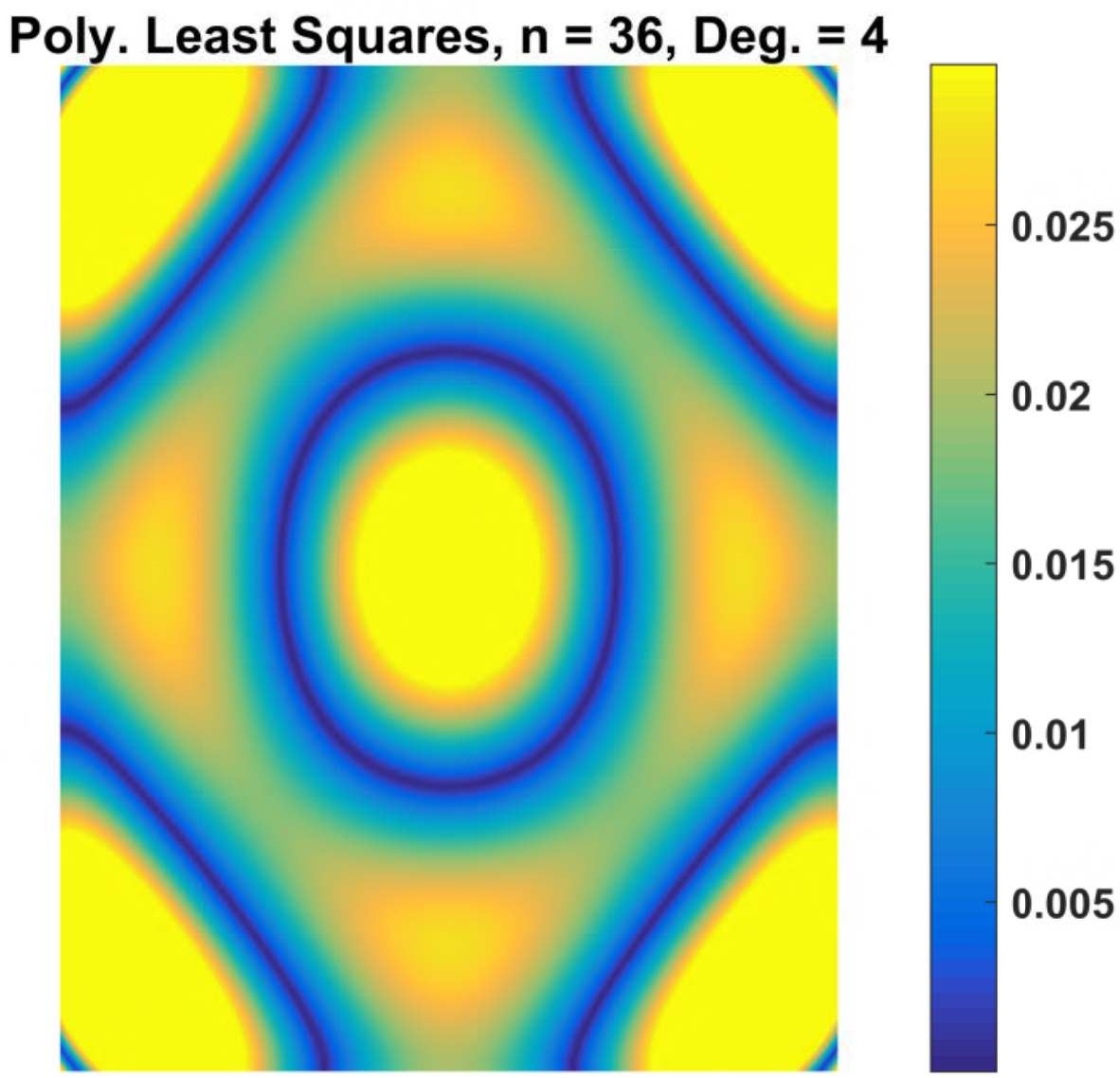} 	
	\label{fig:1a}
}
\subfloat[]
{
	\includegraphics[scale=0.4]{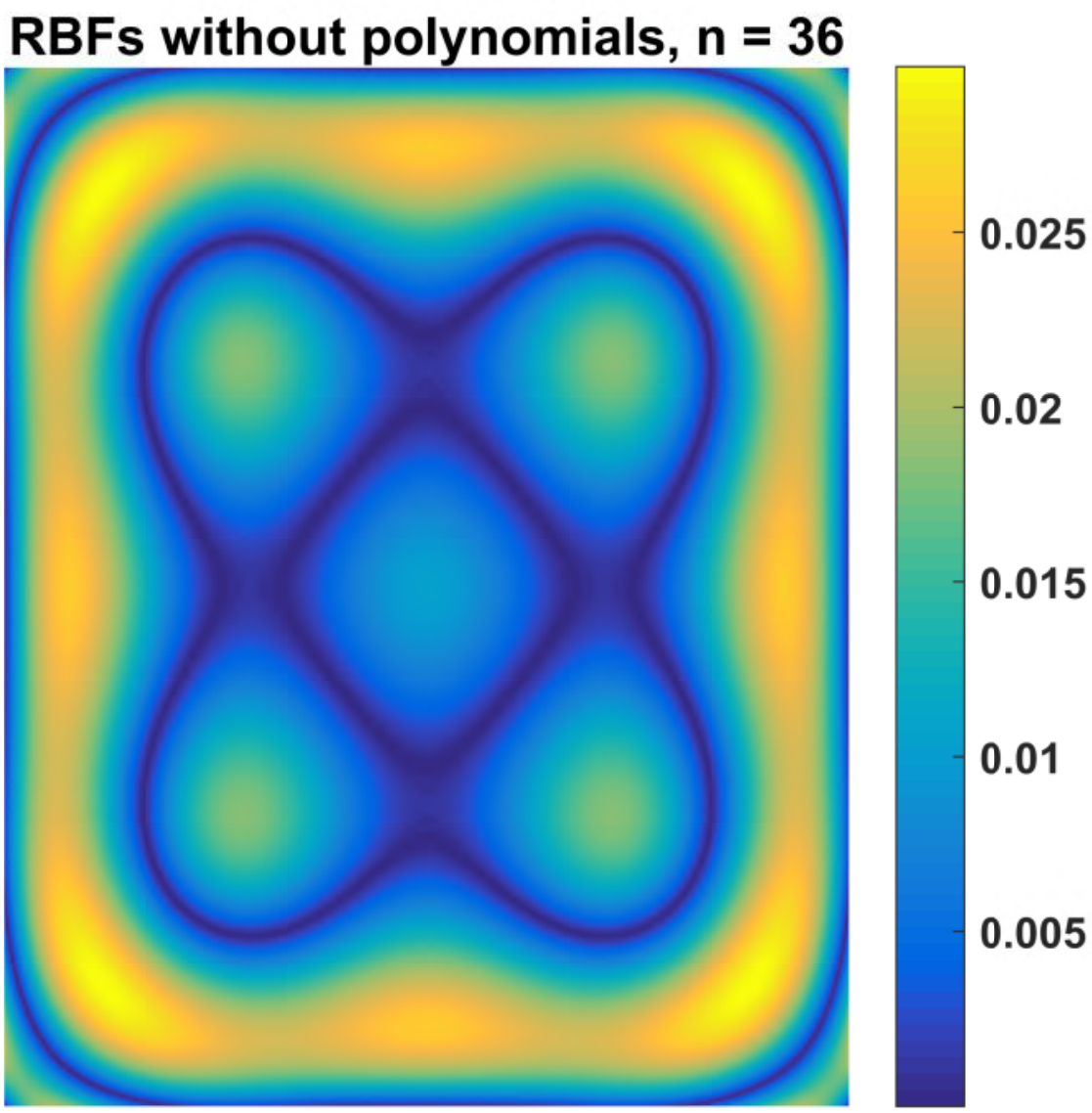}
	
	\label{fig:1b}
}
\subfloat[]
{
	\includegraphics[scale=0.4]{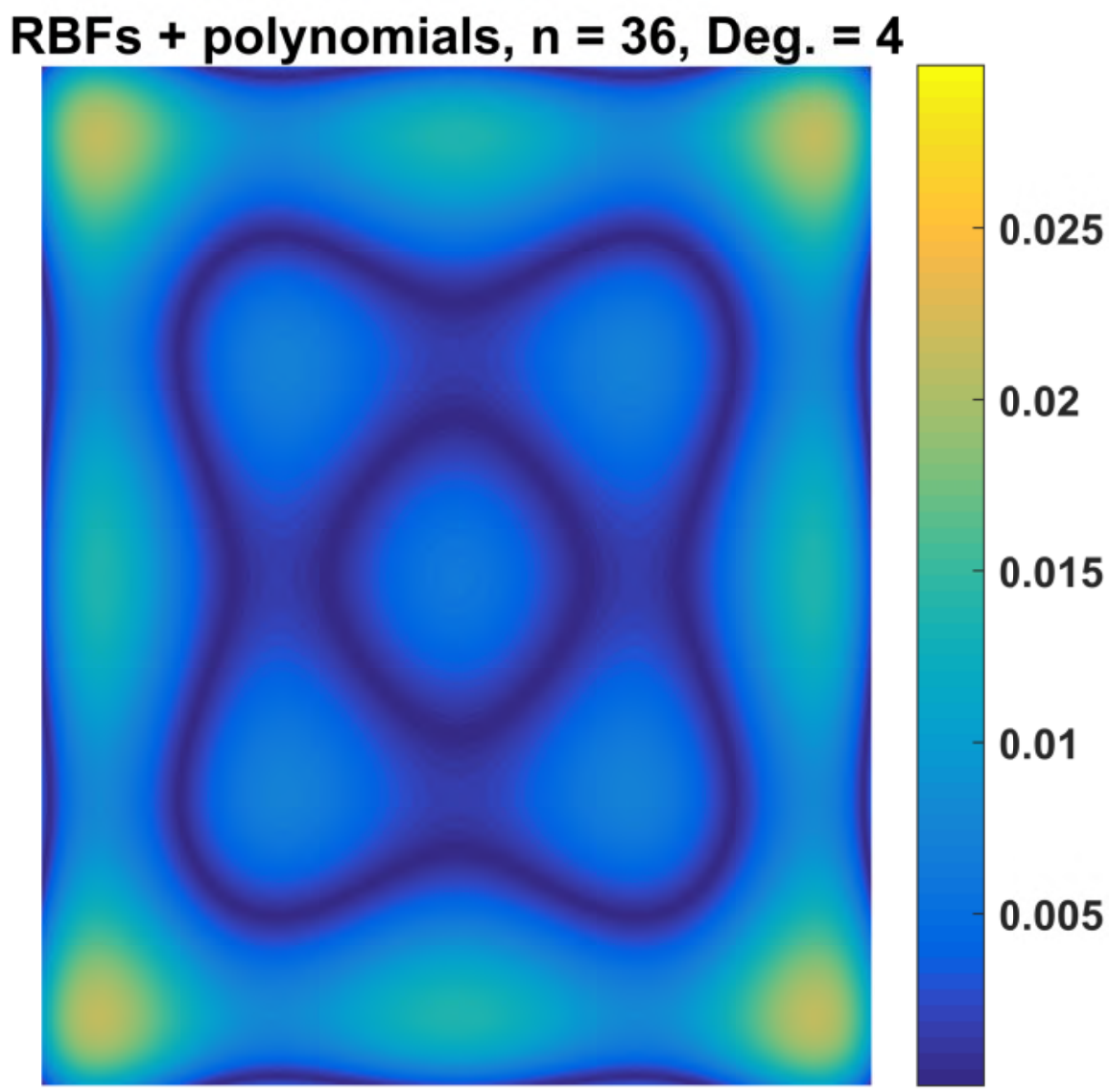} 	
	\label{fig:1c}
}

\subfloat[]
{
	\includegraphics[scale=0.4]{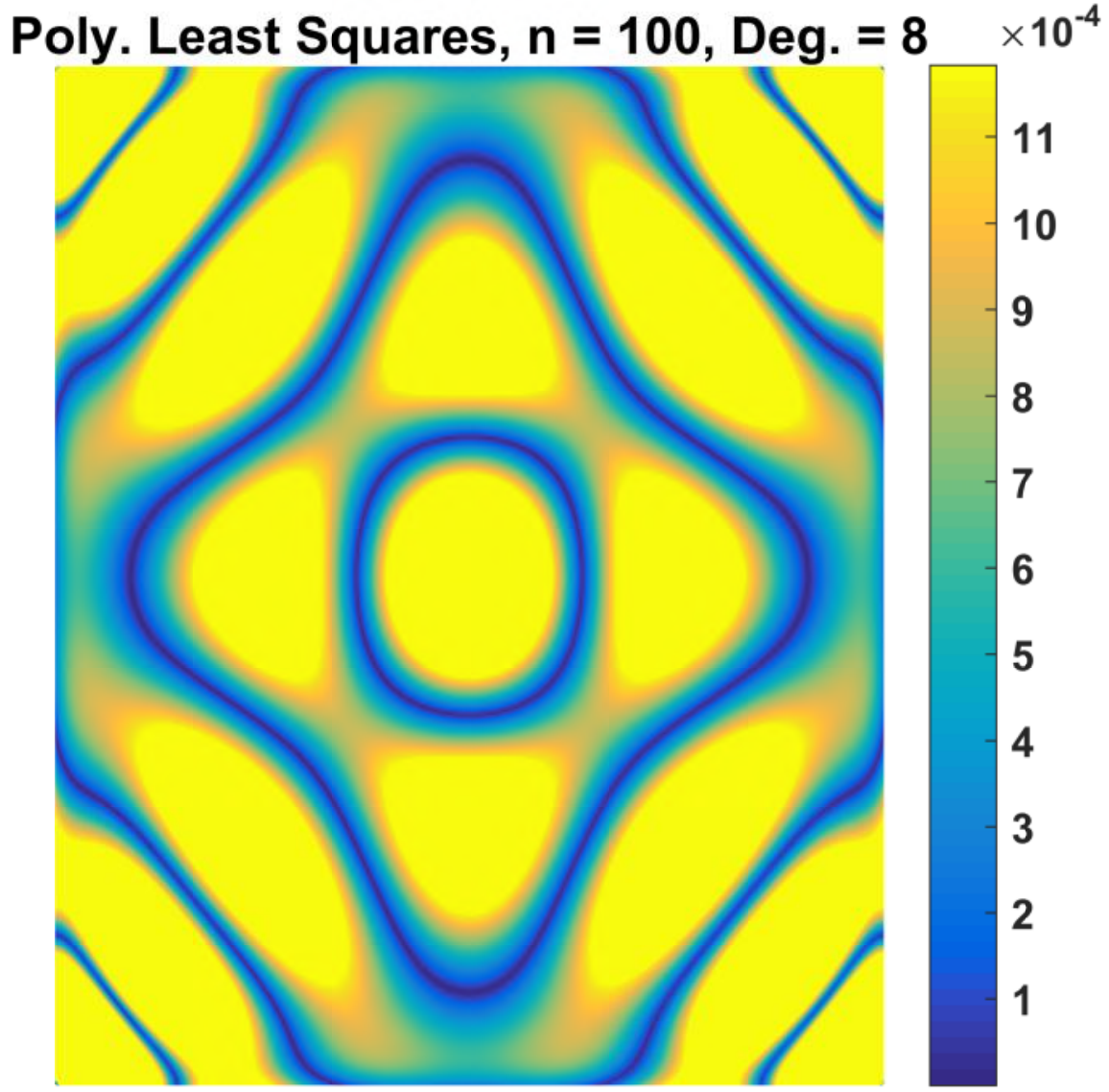} 	
	\label{fig:1g}
}
\subfloat[]
{
	\includegraphics[scale=0.4]{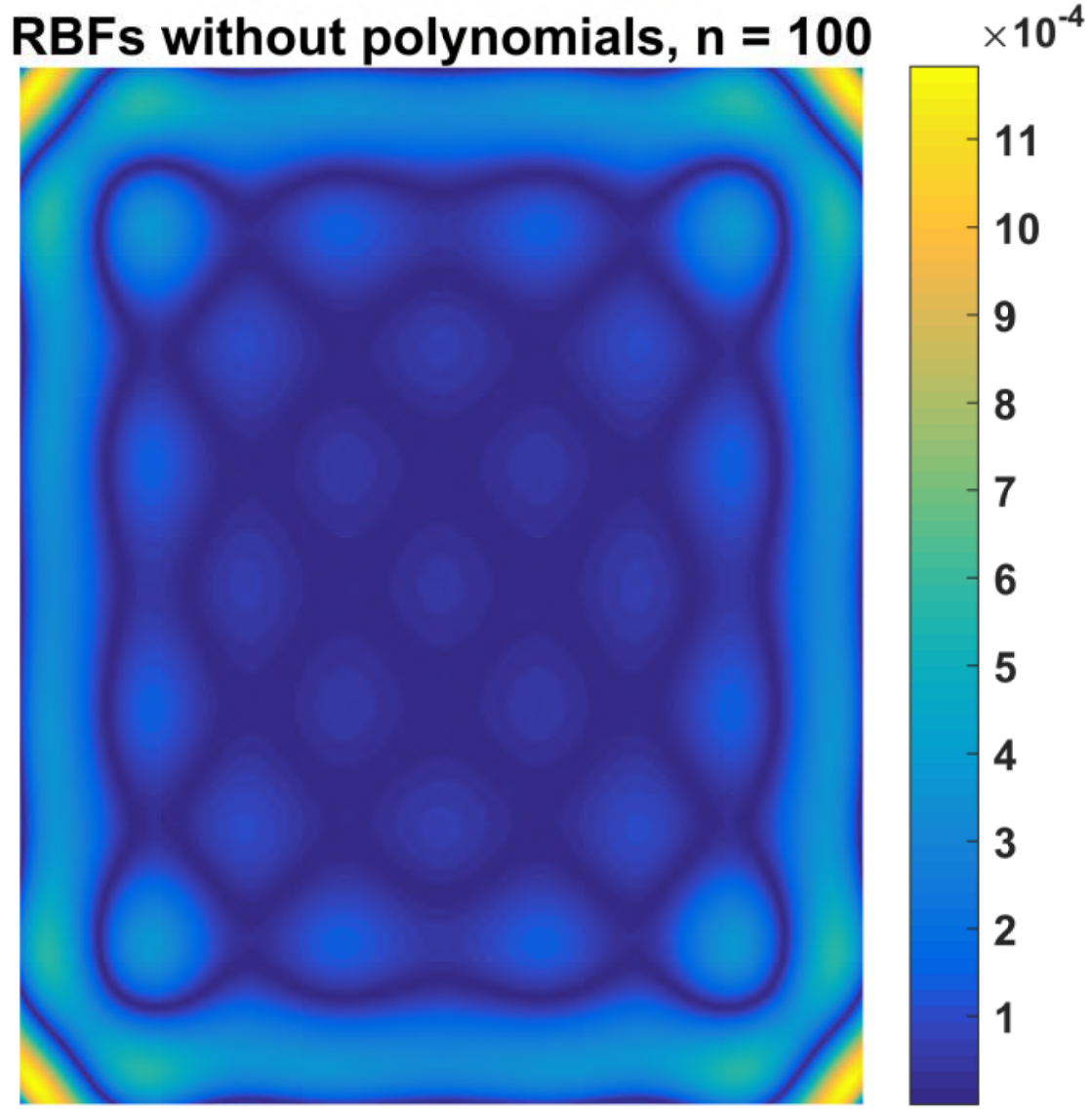}
	
	\label{fig:1h}
}
\subfloat[]
{
	\includegraphics[scale=0.4]{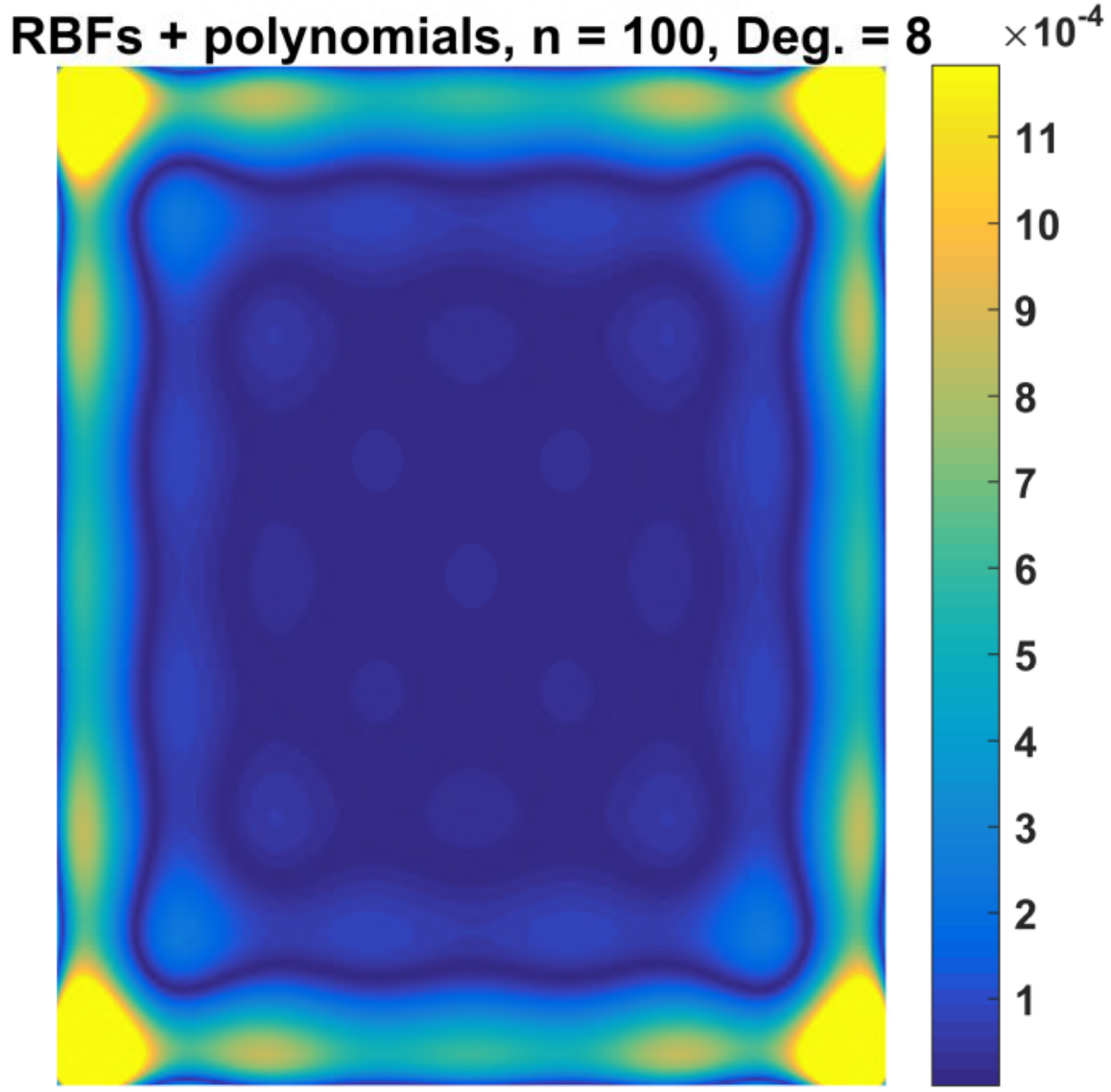} 	
	\label{fig:1i}
}
\caption{Errors in global approximation at Cartesian nodes. The columns show errors in approximating the Runge function with polynomial least squares (left), RBF interpolation (center), and augmented RBF interpolation (right). The rows show errors for increasing numbers of collocation nodes. Darker colors indicate lower errors.}
\label{fig:interp_test1}	
\end{figure}
When interpolating the Runge function using Cartesian nodes, polynomial least squares is known to be useful in overcoming the Runge phenomenon (albeit by sacrificing geometric convergence)~\cite{PlatteImpossible}. On the other hand, global interpolation with either RBFs or augmented RBFs at Cartesian nodes is expected to result in the Runge phenomenon. We interpolate the Runge function at $n=36$ and $n=100$ Cartesian nodes. The results of this test are shown in Figure \ref{fig:interp_test1}.

Evaluation errors decrease in magnitude going from the top to the bottom rows of Figure \ref{fig:interp_test1}. However, the error distribution varies across columns. Figures \ref{fig:1a} and \ref{fig:1g} show that the error is relatively high at the \emph{center} of the domain when using polynomial least squares, in addition to being high at the four corners due to Runge oscillations. It is also easy to see darker zones free of this high error, especially when the polynomial degree is increased. On the other hand, Figures \ref{fig:1b} and \ref{fig:1h} show that the errors for the RBF interpolant are relatively high mainly at the boundaries of the domain. Figure \ref{fig:1b} shows a relatively small dark region in the center of the domain that grows in area as $n$ is increased. Increasing $n$ confines large-amplitude oscillations to the edges of the domain. 

Augmenting the RBFs with the same polynomial degree as in polynomial least squares results in similar error distributions. This can be seen in Figures \ref{fig:1c} and \ref{fig:1i}, with the former showing the greatest improvement. Appending a high degree polynomial to the RBF appears to darken the interior regions of the domain (reducing error there) and lighten the boundary regions (increasing the error there);~\emph{e.g.}, augmenting an interpolant with $n=100$ with polynomials results in Figure \ref{fig:1i}, which has a darker interior region than Figure \ref{fig:1h} (RBFs without polynomials).

This experiment indicates that that RBF interpolants (augmented or otherwise) localize the errors more towards the boundaries of the domain than polynomial least squares. As $n$ is increased, RBF interpolation (augmented or otherwise) confines Runge oscillations to a zone of decreasing size adjacent to the boundaries. Viewing the errors in global RBF interpolation as a \emph{worst-case} stand-in for \emph{per-stencil} RBF-FD errors, this seems to imply that RBF-FD weights likely possess good approximation properties in some finite radius around the centers of each stencil, especially for large $n$. This is contrary to current practice, where RBF-FD weights are computed only at the centers of stencils~\cite{FornbergWrightRunge,FlyerPHS}.

\changer{
\subsection{Halton nodes}
\label{sec:runge_halton}
\begin{figure}[h]
\centering
\subfloat[]
{
	\includegraphics[scale=0.38]{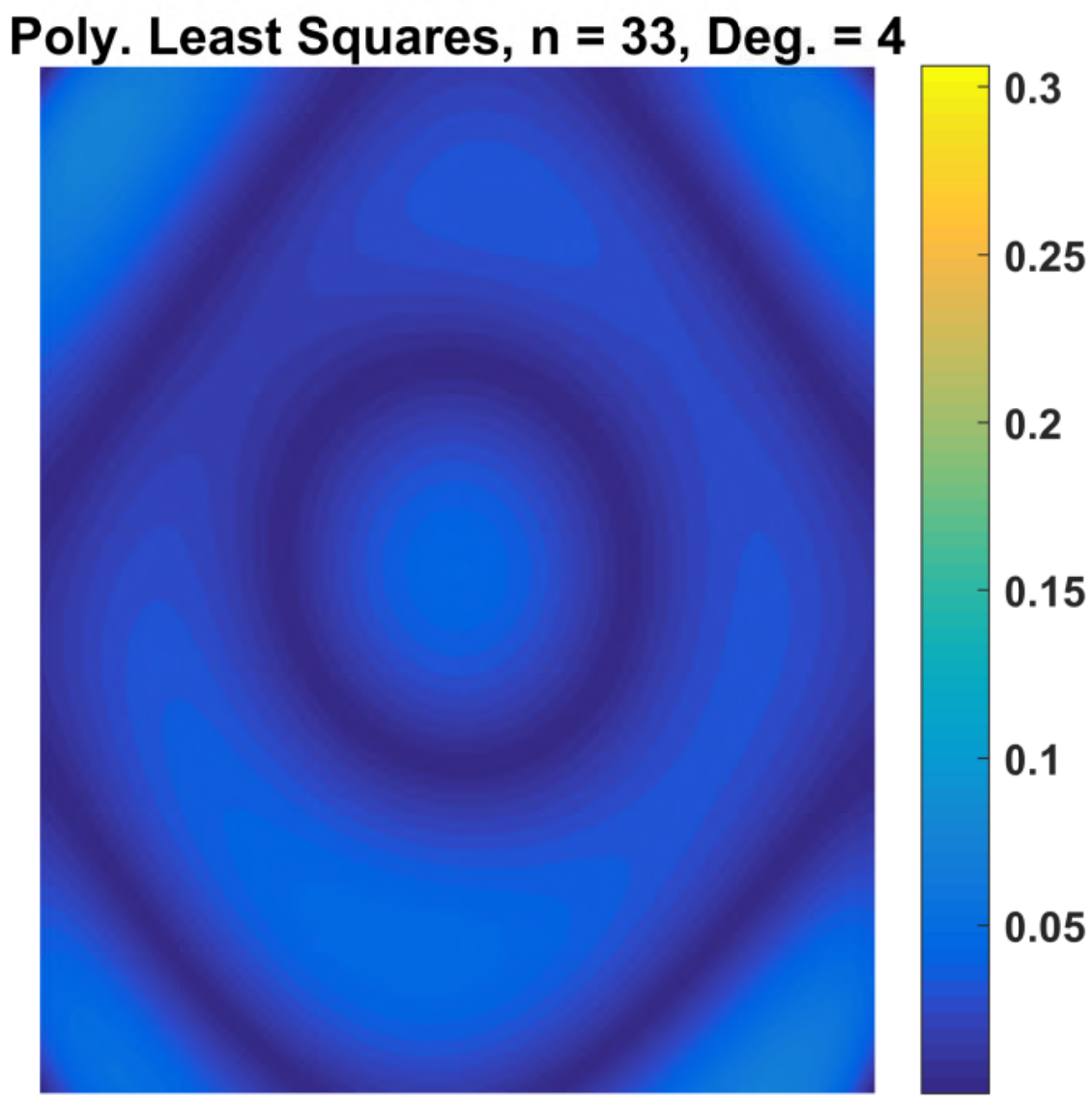} 	
	\label{fig:2a}
}
\subfloat[]
{
	\includegraphics[scale=0.38]{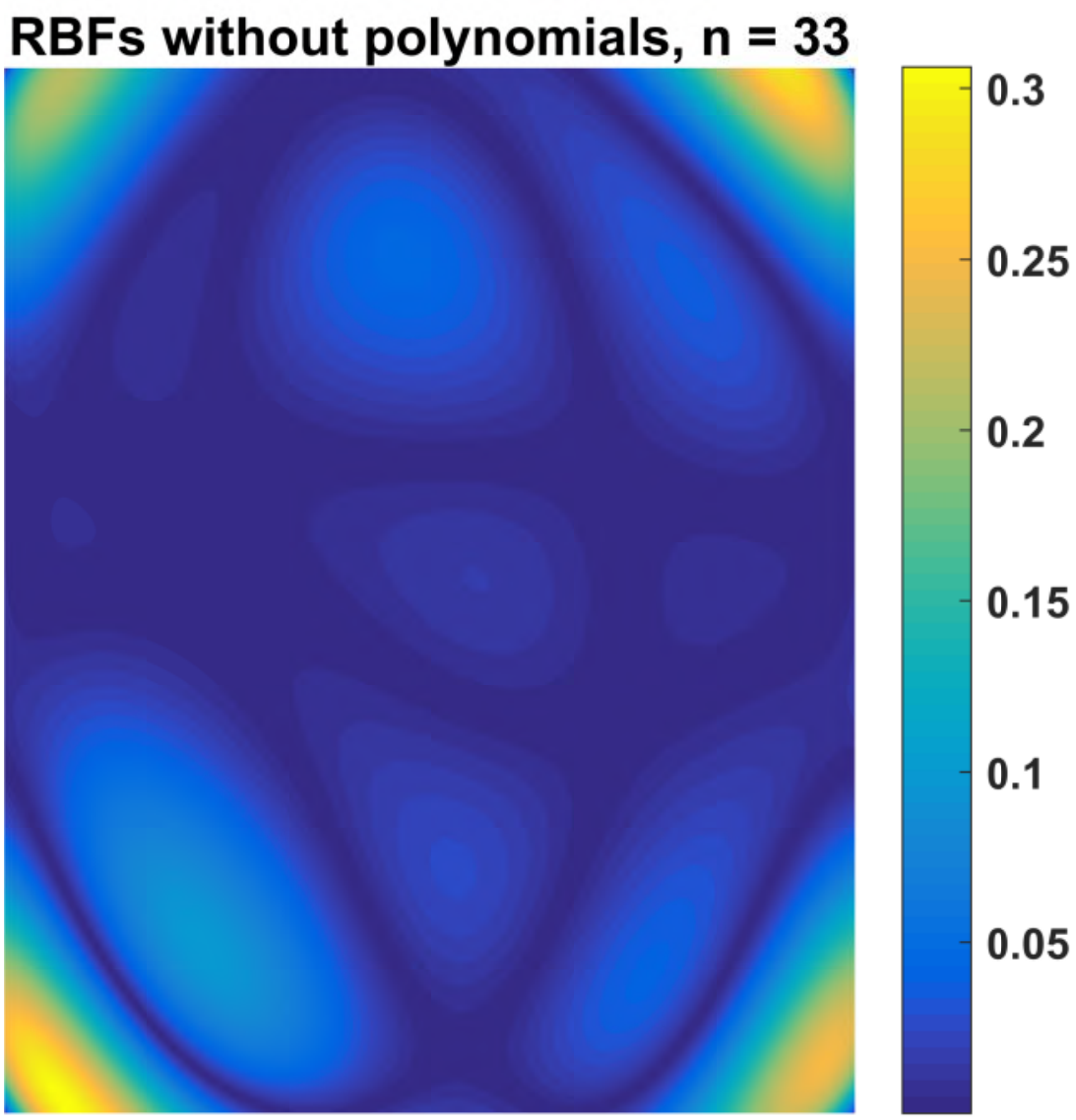}
	
	\label{fig:2b}
}
\subfloat[]
{
	\includegraphics[scale=0.38]{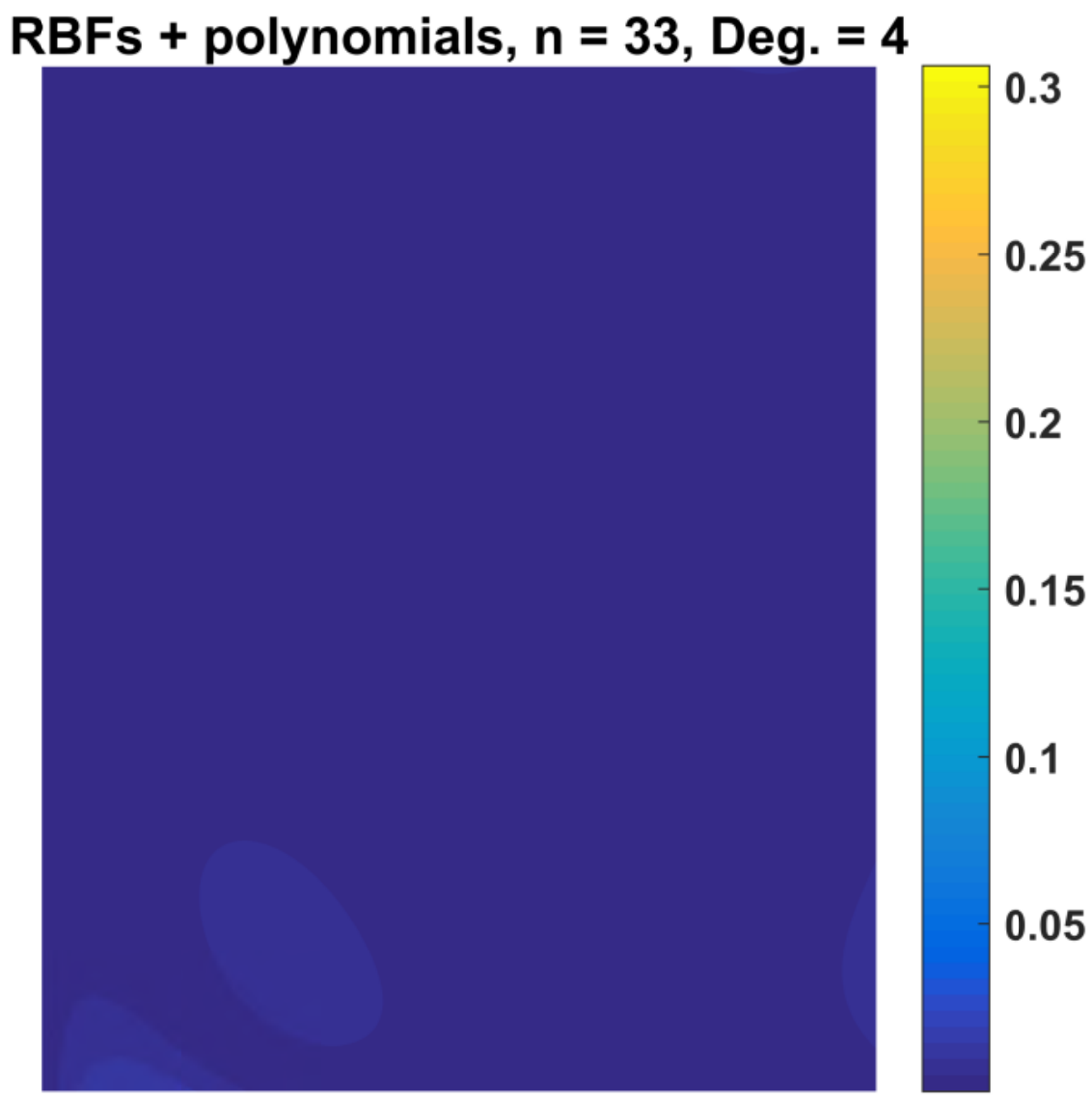} 	
	\label{fig:2c}
}

\subfloat[]
{
	\includegraphics[scale=0.38]{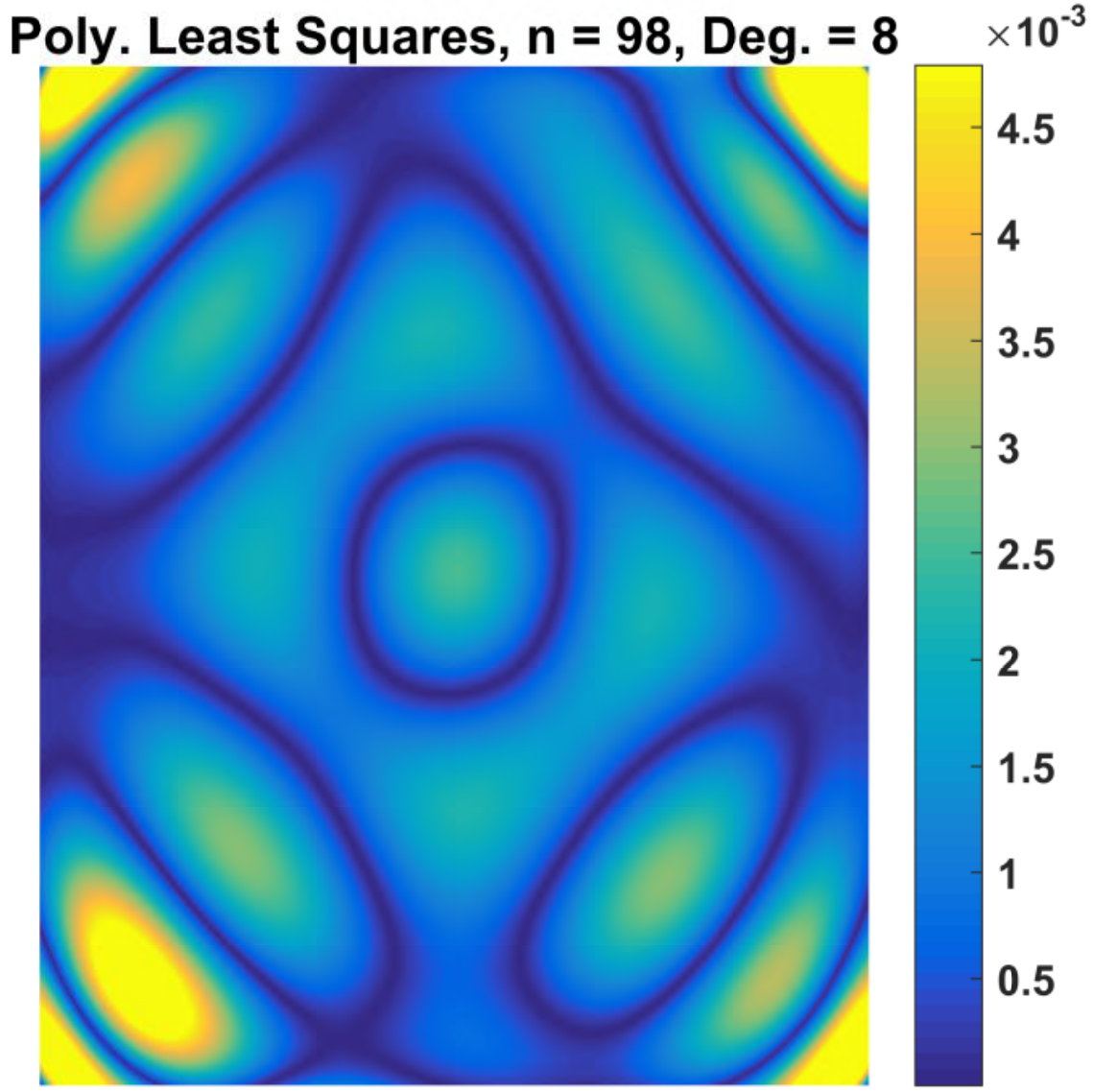} 	
	\label{fig:2g}
}
\subfloat[]
{
	\includegraphics[scale=0.38]{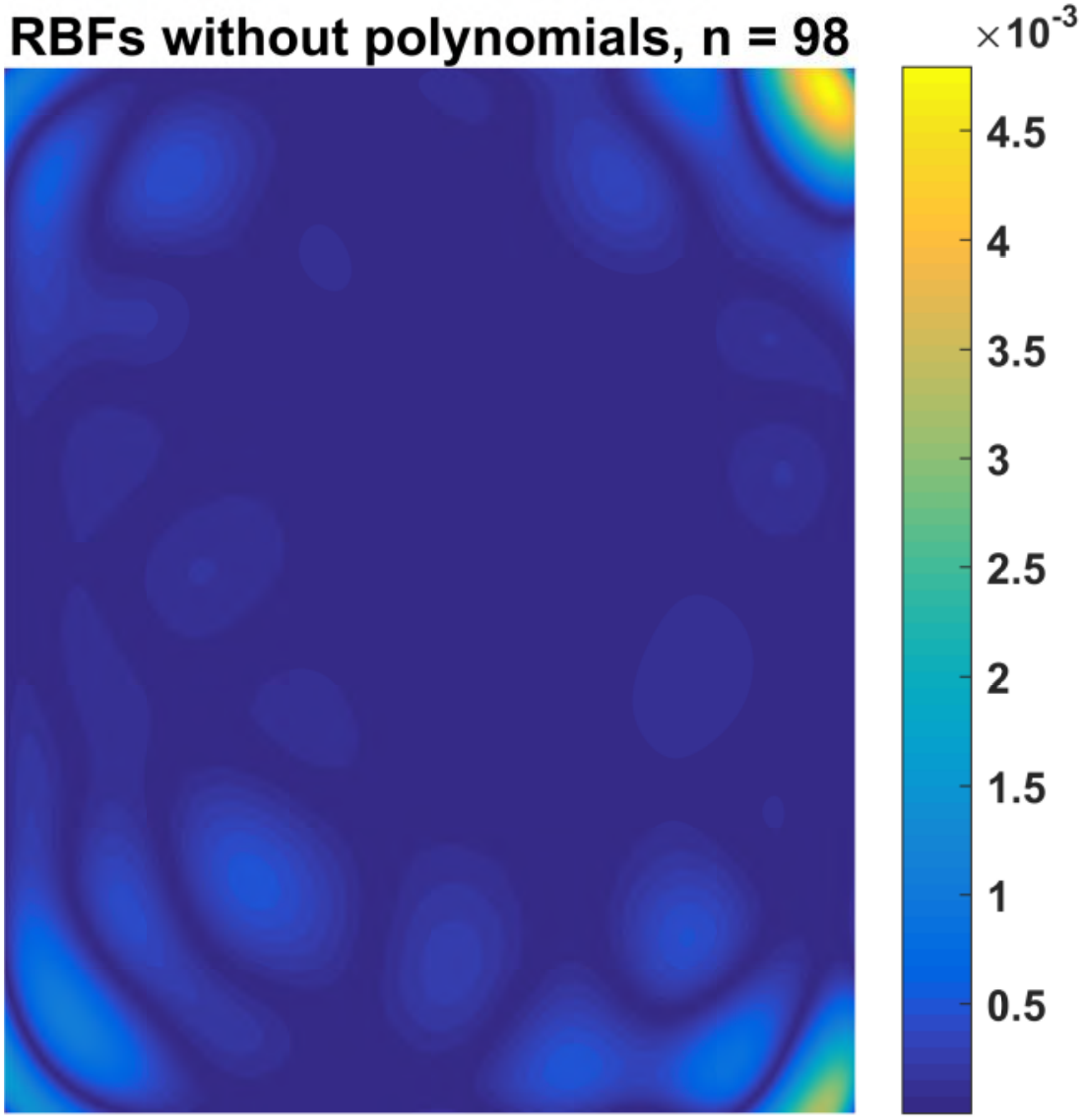}
	
	\label{fig:2h}
}
\subfloat[]
{
	\includegraphics[scale=0.38]{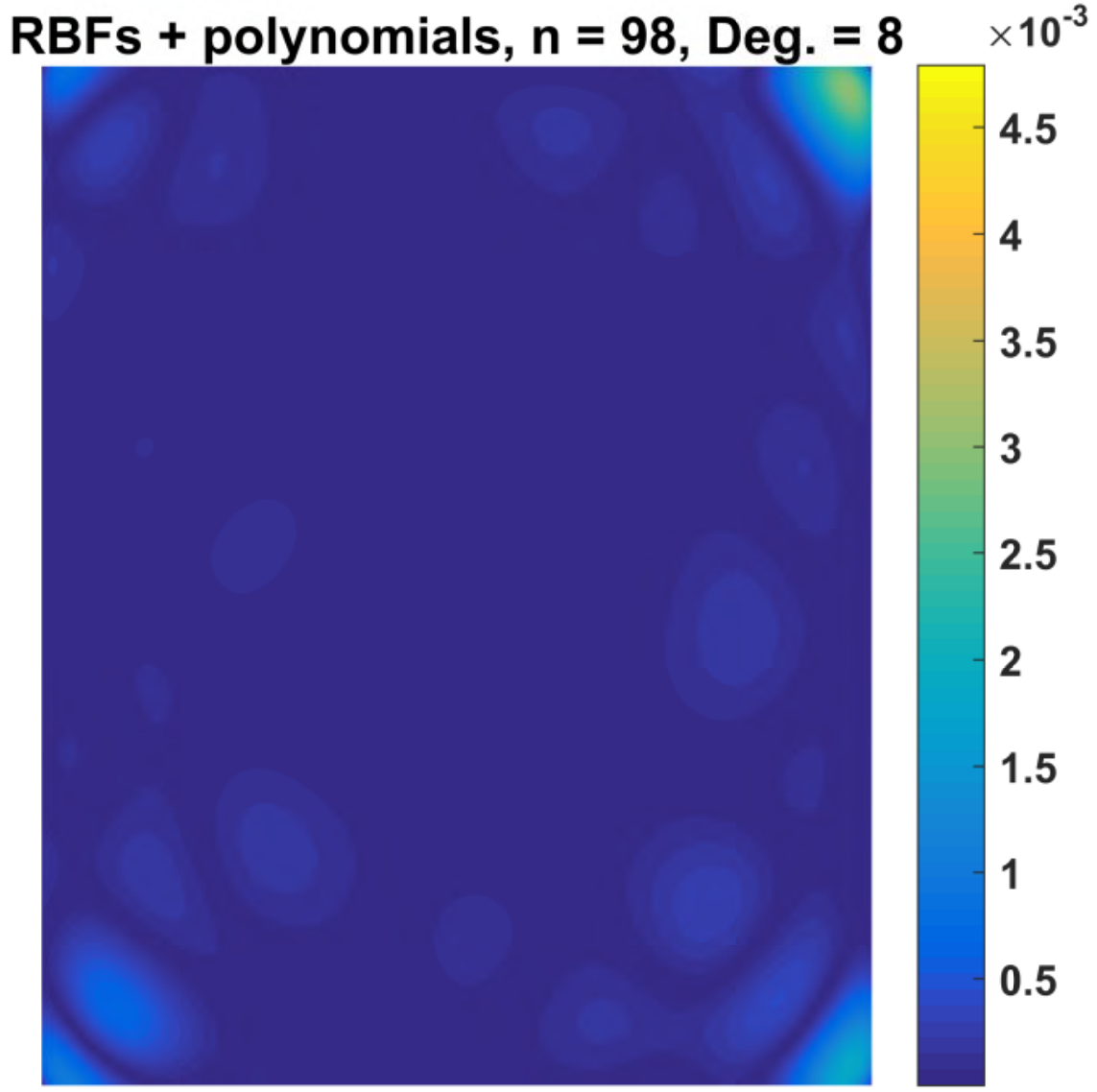} 	
	\label{fig:2i}
}
\caption{Errors in global approximation at Halton nodes. The columns show errors in approximating the Runge function with polynomial least squares (left), RBF interpolation (center), and augmented RBF interpolation (right). The rows show errors for increasing numbers of collocation nodes. Darker colors indicate lower errors.}
\label{fig:interp_test2}	
\end{figure}
We repeat the above experiment using Halton nodes in $(-1,1)^2$ and evenly-spaced nodes on the boundary of $[-1,1]^2$. The results for $n=33$ and $n=98$ are shown in Figure \ref{fig:interp_test2}.

The patterns here are not quite as clear as in the Cartesian case, but it is nevertheless possible to observe some trends. In general, it appears that using Halton nodes localizes errors more towards the boundaries. When $n=33$, we see from Figure \ref{fig:2a} that polynomial least squares still retains its error pattern from Cartesian nodes, but with much lower errors in the interior. The RBF error pattern for $n=33$, on the other hand, does not resemble the Cartesian case (Figure \ref{fig:2b}). The boundary errors are larger for RBF interpolation than for polynomial least squares, but there are larger zones of low error in the interior. Augmented RBF interpolation appears to produce the best error pattern of all, with uniformly low errors up to and including the boundary (Figure \ref{fig:2c}).

Interestingly, for $n=98$, Figure \ref{fig:2g} shows that the errors for polynomial least squares are no longer quite as low in the interior relative to the boundary (compared to $n=33$). In contrast, Figures \ref{fig:2h} and \ref{fig:2i} show that the errors are uniformly low in the interior for RBF and augmented RBF interpolation. For both values of $n$, appending polynomials actually shrinks the boundary error in this scenario.}

\changeb{This experiment shows that node placement plays a vital role in controlling error distribution. For more on the effect of node placement on RBF approximations, we refer the reader to~\cite{FlyerNS}. For the purposes of this article, it is sufficient to note that the errors are low in some region around the center for both Cartesian and Halton nodes. It is possible that these low-error regions correspond to specific features of the function. A thorough exploration of such low-error contours will likely require a Morse-Smale complex~\cite{PascucciMS}, which is built specifically using function features (such as saddles, sources, sinks etc.). We leave an exploration of this to future work.}

\changer{\section{Overlapped RBF-FD}}
\label{sec:overlap}
\begin{figure}[tbhp]
\centering
\includegraphics[scale=0.4]{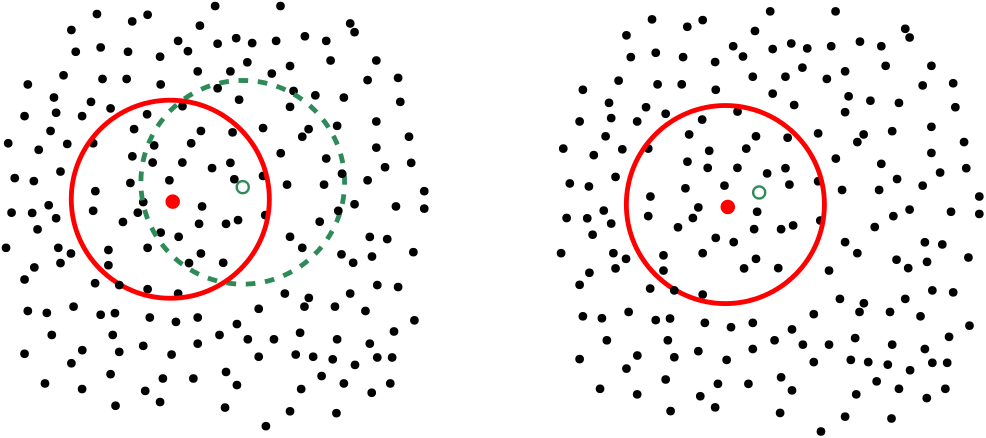} 	
\caption{Illustration of the effect of the overlap parameter $\delta$ in the overlapped RBF-FD method. In the standard RBF-FD method ($\delta = 1$), RBF-FD weights are only computed at the center of each stencil (filled circle and dashed empty circle) and discarded elsewhere (see left). In the overlapped RBF-FD method, depending on $\delta$, the RBF-FD weights for both the filled circle and the dashed empty circle may be computed by the same thick-line stencil (see right).}
\label{fig:orbffd}	
\end{figure}
With evidence supporting the accuracy of global (augmented) RBF interpolation in a region surrounding the center of the domain, we generate the following hypothesis: \emph{RBF-FD weights generated from local RBF interpolation may be of reasonable accuracy in a region around the center of each RBF-FD stencil}. Leaving the testing of this hypothesis to later sections, we now formulate a generalization of the RBF-FD method which we call the \emph{overlapped RBF-FD method}. The terminology ``overlapped'' refers to the extent to which stencils are overlapped; in the standard RBF-FD method, neighboring stencils are considered to be \emph{fully overlapped} in that all weights other than those at the center of these stencils are discarded. In the new method, we will retain weights in some ball around the center of the stencil.

\subsection{Description}
\label{sec:orbffd_desc}
We first describe the augmented RBF-FD method in greater detail. Let $X = \{\vx_k\}_{k=1}^N$ be a global list of nodes on the domain, and without loss of generality, let $\phi(|\vx - \vy\|) = \|\vx - \vy\|^m$, which is the polyharmonic spline (PHS) RBF of order $m$ (where $m$ is odd). Define the stencil $P_k$ to be the set containing $\vx_k$ and its $n-1$ nearest neighbors. As in Section 2, let $\calI^k_1 \hdots \calI^k_n$ be the \emph{global} indices of the nodes in $P_k$, with $k = \calI^k_1$. Now consider the task of computing RBF-FD weights for approximating a linear differential operator $\calL$ on this stencil at the point $\vx_k$ in the stencil.  Following~\cite{FlyerPHS,FlyerNS}, this can be written as:
{\footnotesize
\begin{align}
\underbrace{
\begin{bmatrix}
\|\vx_{\calI^k_1} - \vx_{\calI^k_1}\|^m & \hdots & \|\vx_{\calI^k_1} - \vx_{\calI^k_n}\|^m & \psi^k_1(\vx_{\calI^k_1}) & \hdots & \psi^k_M(\vx_{\calI^k_1}) \\
\vdots & \ddots & \vdots & \vdots & \ddots & \vdots \\
\|\vx_{\calI^k_n} - \vx_{\calI^k_1}\|^m & \hdots & \|\vx_{\calI^k_n} - \vx_{\calI^k_n}\|^m & \psi^k_1(\vx_{\calI^k_n}) & \hdots & \psi^k_M(\vx_{\calI^k_n}) \\
\psi^k_1(\vx_{\calI^k_1}) & \hdots & \psi^k_1(\vx_{\calI^k_n}) & 0 &\hdots & 0 \\
\vdots & \ddots & \vdots & \vdots & \ddots & \vdots\\
\psi^k_M(\vx_{\calI^k_1}) & \hdots & \psi^k_M(\vx_{\calI^k_n}) & 0 & \hdots & 0
\end{bmatrix}}_{\hat{A}_k}
\underbrace{
\begin{bmatrix}
w^{\phi}_{1} \\
\vdots \\
w^{\phi}_{n} \\
w^{\psi}_1 \\
\vdots\\
w^{\psi}_M
\end{bmatrix}}_{W_k}
= 
\underbrace{
\begin{bmatrix}
\lf.\calL\|\vx - \vx_{\calI^k_1}\|^m\rt|_{\vx = \vx_k}\\
\vdots\\
\lf.\calL\|\vx - \vx_{\calI^k_n}\|^m\rt|_{\vx = \vx_k}\\
\lf.\calL \psi^k_1(\vx)\rt|_{\vx = \vx_k}\\
\vdots\\
 \lf.\calL\psi^k_M(\vx)\rt|_{\vx = \vx_k}
\end{bmatrix}}_{B_k}.
\label{eq:rbffd_linsys}
\end{align}
}
\changer{The upper left block of $\hat{A}_k$ is $A_k$ from Section 2, the upper right block is $\Psi_k$}. Here $\psi_j(\vx)$ is the $j$-th monomial corresponding to the multivariate total degree polynomial in $d$ dimensions. Partition $W_{k}$ as
\begin{align}
W_{k} = \lf[W^{\phi,k},W^{\psi,k}\rt]^T.
\end{align}
Then, $W^{\phi,k}$ is an $n \times 1$ matrix of RBF-FD weights for the corresponding node $j$ in the stencil. $W^{\psi,k}$ can be safely discarded~\cite{FlyerPHS}.

The overlapped RBF-FD method adds columns to the right hand side matrix $B_k$. Before doing so, we define the \emph{stencil width} $\rho_k$ as
\begin{align}
\rho_k = \max\limits_{j} \|\vx_k - \vx_{\calI^k_j}\|, j=1,\hdots,n.
\end{align}
Given a parameter $\delta \in [0,1]$, this allows us to define the \emph{stencil retention distance} $r_k$ as
\begin{align}
r_k = (1-\delta)\rho_k.
\end{align}
Here, $\delta$ is called the \emph{overlap parameter}. The parameters $\rho_k$ and $\delta$ thus effectively define a \emph{stencil retention ball} $\mathbb{B}_k$ of radius $r_k$ centered at each node $\vx_k$. When $\delta=1$, the ball $\mathbb{B}_k$ collapses to a single point--its center $\vx_k$.  Now, let $R_k$ be the set of \emph{global indices} of the nodes in the subset $\mathbb{B}_k \subseteq P_k$:
\begin{align}
R_k = \{\calR^k_1, \calR^k_2, \hdots, \calR^k_{p_k}\},
\end{align}
where $1 \leq p_k \leq n$. In general, $R_k$ is some permutation of a \emph{subset} of the indices of the nodes in $P_k$. There are two important cases: a) if $\delta = 1$, $R_k = \calI^k_1$ and only $\vx_{\calI^k_1}$ lies in $\mathbb{B}_k$ (standard RBF-FD method); and b) if $\delta = 0$, all the nodes $\vx_{\calI^k_j}$,$j=1,\hdots,n$ lie in $\mathbb{B}_k$.

The overlapped RBF-FD method modifies both $B_k$ and $W_k$ by adding columns corresponding to all the nodes within $\mathbb{B}_k$ so that we compute RBF-FD weights at all nodes that lie within the retention distance $r_k$. We define new matrices $\hat{W}_k$ and $\hat{B}_k$ so that
\begin{align}
\hat{A}_k \hat{W}_k = \hat{B}_k,
\label{eq:orbffd_eq}
\end{align}
where
\begin{align}
\hat{B}_k = 
\begin{bmatrix}
\lf.\calL\|\vx - \vx_{\calI^k_1}\|^m\rt|_{\vx = \vx_{\calR^k_1}} & \hdots & \lf.\calL\|\vx - \vx_{\calI^k_1}\|^m\rt|_{\vx = \vx_{\calR^k_{p_k}}}\\
\vdots & \ddots & \vdots \\
\lf.\calL\|\vx - \vx_{\calI^k_n}\|^m\rt|_{\vx = \vx_{\calR^k_1}} & \hdots & \lf.\calL\phi\|\vx - \vx_{\calI^k_n}\|^m\rt|_{\vx = \vx_{\calR^k_{p_k}}} \\
\lf.\calL \psi^k_1(\vx)\rt|_{\vx = \vx_{\calR^k_1}} & \hdots & \lf.\calL \psi^k_1(\vx)\rt|_{\vx = \vx_{\calR^k_{p_k}}} \\
\vdots & \ddots & \vdots \\
 \lf.\calL \psi^k_M(\vx) \rt|_{\vx = \vx_{\calR^k_1}} & \hdots & \lf.\calL \psi^k_M(\vx) \rt|_{\vx = \vx_{\calR^k_{p_k}}}
\end{bmatrix}.
\label{eq:orbffd_rhs}
\end{align}
$\hat{B}_k$ and $\hat{W}_k$ are now $(n + M) \times p_k$ matrices, where $1 \leq {p_k} \leq n$. As before, partition $\hat{W}_k$ as $\hat{W}_k = \lf[\hat{W}^{\phi,k}, \hat{W}^{\psi,k}\rt]^T$: the $n \times {p_k}$ matrix $\hat{W}^{\phi,k}$ now contains the $n$ RBF-FD weights for each of the ${p_k}$ nodes in the retention ball, and $\hat{W}^{\psi,k}$ is again discarded. If $\delta = 1$, we have $\hat{B}_k = B_k$, $\hat{W}_k = W_k$, and we recover the RBF-FD method. On other hand, if $\delta = 0$, $\hat{W}_k$ contains weights for every point in the stencil. In this work, we never allow the value of $\delta = 0$ as this always appears to result in instabilities in the context of PDEs. 

As in the augmented RBF-FD method, the local RBF-FD weights are stored in the rows of a global \emph{sparse} differentiation matrix $L$ using the stencil index sets $\calI^k$. More specifically, we have
\begin{align}
L_{\calR^k_i, \calI^k_j} = \hat{W}^{\phi,k}_{j,i}, i=1,\hdots,p_k, j=1,\hdots,n.
\end{align}
As presented thus far, the overlapped RBF-FD method generates multiple candidate sets of entries for the rows in $L$ that are shared across stencils. To only generate a single set of entries per row of $L$, we simply require that weights computed for any node $\vx_k$ never be recomputed again by some other stencil $P_i$, $i\neq k$. The immediate consequence of this approach is that the order in which the nodes are traversed determines the RBF-FD weights assigned to a node and its stencil constituents. Our experiments indicate that this is not detrimental to the solution of PDEs. A second consequence of the requirement not to repeat weight computations is that the overlapped RBF-FD method uses fewer than $N$ stencils for $N$ global nodes. In other words, if $N_{\delta}$ is the number of stencils in the overlapped method, $N_{\delta} < N$. 

\subsection{Local Lebesgue functions for error estimation and improved stabilty}
\label{sec:leb}

Up to this point, we have made no mention of the formal errors or stability of our technique. Based on Figure 1, one could anticipate that very small values of $\delta$ could cause $\mathbb{B}_k$ to encompass nodes in the ``Runge zone''; in the context of RBF-FD stencils, this means that we may compute weights that result in both high errors and numerical instability. This instability manifests as spurious eigenvalues in the spectrum of the differentiation matrix. For this article, this means that when approximating the Laplacian $\calL \equiv \Delta$, the differentiation matrix $L$ may have eigenvalues with positive real parts. The goal of this section is twofold: first, to make the overlapped RBF-FD method more robust for small values of $\delta$ and relate the eigenvalues of $L$ to the RBF-FD weights; and second, to discuss the errors associated with the augmented RBF-FD method (overlapped or otherwise).  

\subsubsection{Improving stability}
\label{sec:stab_leb}
While it is not practical to require stability for all values of $\delta$, we will nevertheless present an approach that appears to improve stability for $\delta \geq 0.2$. In general, we caution against selecting $\delta < 0.2$.

To aid the discussion, consider the \emph{Lagrange} form of the augmented RBF interpolant. Given a function $f$ that we wish to approximate on the node set $X$, we now write the augmented RBF interpolant on the stencil $P_k$ as:
\begin{align}
s^k(\vx) = \sum\limits_{j=1}^n \ell^k_j(\vx) f_{\calI^k_j},
\end{align}
where $\ell^k_j (\vx)$ are the \emph{local Lagrange functions} or \emph{cardinal functions} on the stencil $P_k$, and $f_{\calI^k_j}$ are samples of $f$ on the stencil $P_k$ at the nodes $\vx_{\calI^k_j}$. The cardinal functions have the Kronecker delta property: $\ell^k_j(\vx_{\calI^k_j}) = 1$, $\ell^k_j(\vx_{\calI^k_i}) = 0, \forall i \neq j$. The derivatives of the local Lagrange functions give the RBF-FD weights on this stencil~\cite{Wright200699}, and their integrals can be used to generate quadrature rules~\cite{LLF}. Unlike in~\cite{LLF}, our goal is not to use the local Lagrange functions as approximants, but rather to use them to develop a stabilization procedure.

Without explicitly computing $\ell^k_j$, we can use it to define the \emph{local Lebesgue function} $\Lambda^k(\vx)$ for the stencil $P_k$ as
\begin{align}
\Lambda^k(\vx) = \sum\limits_{j=1}^n | \ell^k_j (\vx)|.
\end{align}
We use an analogue of the local Lebesgue function to help develop a stability indicator and stabilization procedure for the overlapped RBF-FD method. Consider applying a linear differential operator $\calL$ to the Lagrange form of the augmented RBF interpolant on the stencil $P_k$:
\begin{align}
\calL s^k(\vx) = \sum\limits_{j=1}^n \lf(\calL \ell^k_j(\vx) \rt) f_{\calI^k_j}.
\end{align}
Written in this form, the quantities $\calL \ell^k_j(\vx)$ are the entries of the matrix $\hat{W}^{\phi}_k$ arising from partition $\hat{W}_k$ in~\eqref{eq:orbffd_eq}. We can now define an analogue to the Lebesgue function that corresponds to the differential operator $\calL$: the local $\calL$-Lebesgue function, denoted by $\Lambda^k_{\calL}$. This function is given explicitly by
\begin{align}
\Lambda^k_{\calL}(\vx) = \sum\limits_{j=1}^n |\calL \ell^k_j(\vx)|.
\end{align}
We evaluate $\Lambda^k_{\calL}$ on the stencil $P_k = \{ \vx_{\calI^k_j} \}_{j=1}^n$ to obtain a pointwise stability indicator corresponding to the operator $\calL$. More specifically, we compute RBF-FD weights for $\calL$ at the point $\vx_{\calI^k_j}$ using the stencil $P_k$ only if
\begin{align}
\Lambda^k_{\calL}(\vx_{\calI^k_j}) \leq \Lambda^k_{\calL}(\vx_k),
\label{eq:calLstab}
\end{align}
where $\vx_k$ is the node closest to the centroid of $P_k$. If this condition is not satisfied, the RBF-FD weights for the node $\vx_{\calI^k_j}$ will be computed using a different stencil $P_i$, $i \neq k$. While we considered other possibilities for stability indicators such as the native space norm of the interpolant~\cite{iske2002}, this indicator appears to be the most robust in the context of the overlapped method.

We now explain the connection between $\Lambda^k_{\calL}$ and the eigenvalues of the differentiation matrix $L$. Using the fact that derivatives of the local Lagrange functions are the RBF-FD weights, we have the following relation:
\begin{align}
\Lambda^k_{\calL}(\vx_{\calI^k_j}) = \sum\limits_{i=1}^n |\hat{W}^{\phi,k}_{ij}|.
\end{align}
More simply, define the vector ${\bf w}^{k,j}$ as:
\begin{align}
{\bf w}^{k,j} = \lf[ \hat{W}^{\phi,k}_{1j},\hat{W}^{\phi,k}_{2j},\hdots,\hat{W}^{\phi,k}_{nj} \rt]^T,
\end{align}
the contents of the $j$\textsuperscript{th} column of $\hat{W}^{\phi,k}$. Then,
\begin{align}
\Lambda^k_{\calL}(\vx_{\calI^k_j}) = \| {\bf w}^{k,j} \|_1.
\end{align}
Now, by definition, the rows of the differentiation matrix $L$ contain the RBF-FD weights. Assume without loss of generality that the $k$\textsuperscript{th} row of $L$ contains the vector ${\bf w}^{k,1}$ distributed through its columns, interspersed with zeros. Recall that for the $k$\textsuperscript{th} row, the radius of the $k$\textsuperscript{th} \emph{Gershgorin disk} is given by the sum of the off-diagonal elements of the $k$\textsuperscript{th} row:
\begin{align}
g_k = \sum\limits_{j\neq k} |L_{kj}|, j=1,\hdots,N.
\end{align}
Most of the row entries are zero, with the only non-zero entries being the RBF-FD weights for that row. Assuming without loss of generality that the diagonal element of the $k$\textsuperscript{th} row is given by $|\hat{W}^{\phi,k}_{11}|$, we have
\begin{align}
g_k = \|{\bf w}^{k,1}\|_1 - |\hat{W}^{\phi,k}_{11}|.
\end{align}
This in turn can be written as
\begin{align}
g_k = \Lambda^k_{\calL}(\vx_k) - |\hat{W}^{\phi,k}_{11}|.
\end{align}
\emph{i.e.}, the radius of the $k$\textsuperscript{th} Gershgorin disk depends on the local $\calL$-Lebesgue function. From the Gershgorin circle theorem, we know that the eigenvalues of $L$ are contained within the union of all Gerschgorin disks of $L$. If $\lambda$ is an eigenvalue of $L$, we have
\begin{align}
|\lambda - L_{kk}| &\leq g_k, \\
\implies |\lambda - \hat{W}^{\phi,k}_{11}| &\leq \Lambda^k_{\calL}(\vx_k) - |\hat{W}^{\phi,k}_{11}|.
\end{align}
If $-\Lambda^k_{\calL}(\vx_k) < \hat{W}^{\phi,k}_{11} \leq 0$, we have 
\begin{align}
\lambda \in [-\Lambda^k_{\calL}(\vx_k), \Lambda^k_{\calL}(\vx_k) + 2\hat{W}^{\phi,k}_{11}].
\end{align}
A sufficient condition for eigenvalues with non-positive real parts is then:
\begin{align}
2\hat{W}^{\phi,k}_{11} \geq -\Lambda^k_{\calL}(\vx_k),
\end{align}
assuming that the weights $\hat{W}^{\phi,k}_{11}$ are real. On the other hand, if $0<\hat{W}^{\phi,k}_{11}<\Lambda^k_{\calL}(\vx_k)$, we have
\begin{align}
\lambda \in [2\hat{W}^{\phi,k}_{11} - \Lambda^k_{\calL}(\vx_k), \Lambda^k_{\calL}(\vx_k)],
\end{align}
which is unfortunately not useful in generating a sufficiency condition. In general, the approach given by \eqref{eq:calLstab} works regardless of whether the weights are negative by selecting weights that produce smaller values of $\Lambda^k_{\calL}(\vx_k)$ than if no stabilization had been used. The effect of stabilization is discussed in Section \ref{sec:stab}. In general, $\Lambda^k_{\calL}(\vx)$ depends on the stencil node set and the operator being approximated. In practice, we find that for values of $\delta$ that do not extend the retention balls $\mathbb{B}_k$ into the Runge zones, the overlapped RBF-FD method produces a discrete Laplacian $L$ whose eigenvalues have negative real parts and relatively small non-zero imaginary parts \emph{if augmented RBF-FD does the same}. 

\subsubsection{Error Estimates}
We now discuss error estimates for the augmented RBF-FD method that also apply to the overlapped method. The following discussion summarizes the numerical differentiation error estimates developed by Davydov and Schaback~\cite{DavydovMinimal2016,DavydovSchaback2016}, adapted to our notation. These estimates also involve the local $\calL$-Lebesgue functions.

Let $\Omega$ be the domain where we are approximating the differential operator $\calL$. Further, partition $\Omega$ into $N_{\delta}$ sub-domains $\Omega_k$ so that $\Omega = \bigcup\limits_{k=1}^{N_{\delta}} \Omega_k$, where each $\Omega_k$ is the convex hull of the stencil $P_k$. The Sobolev space $W^{r,p}(\Omega_k)$ is given by:
\begin{align}
W^{r,p}(\Omega_k) = \{ f \in L^p(\Omega_k): D^{\alpha}f \in L^p(\Omega) \ \forall \ |\alpha| \leq r \},
\end{align}
where $\alpha \in \mathbb{Z}^d_+$ is a multi-index. The Sobolev $\infty$ norm is then defined as
\begin{align}
\|f\|_{W^{r,\infty}(\Omega_k)} := \max\limits_{|\alpha|\leq r} \|D^{\alpha} f \|_{L^{\infty}(\Omega_k)}.
\end{align}
Let $C^{r,\gamma}(\Omega_k)$ denote the \Holder \ space consisting of all $r$-times continuously differentiable functions $f$ on $\Omega_k$ such that $\|D^{\alpha} f\|_{\gamma} < \infty$ with $|\alpha| = r$. Here, the seminorm $\|.\|_{\gamma}$ is defined for some function $g$ as
\begin{align}
\|g \|_{\gamma} := \sup_{\vx \neq \vy} \frac{|g(\vx) - g(\vy)|}{\|\vx - \vy \|^{\gamma}_2},
\end{align}
for $\vx,\vy \in \Omega_k$. This is a seminorm on $C^{0,\gamma}(\Omega)_k)$. If $f$ and its derivatives up to order $r$ are bounded on $\Omega_k$, we define the semi-norm
\begin{align}
\| f\|_{C^{r,\gamma}} := \max\limits_{|\alpha| = r} \|D^{\alpha} f\|_{\gamma}.
\end{align}
If $\Omega_k$ has a $C^1$ boundary, we can use Morrey's inequality on the right hand side to convert $\|.\|_{\gamma}$to a more familiar Sobolev (semi)norm:
\begin{align}
\| f\|_{C^{r,\gamma}} \leq T\max\limits_{|\alpha| = r}\|D^{\alpha} f\|_{W^{1,\infty}(\Omega_k)},
\end{align}
where $T$ is some constant. We can now use these definitions to rewrite Eq. (13) from~\cite{DavydovMinimal2016}. Let $\vx_k \in \Omega_k$ be a point in the node set $X$ at which we wish to approximate $\calL f$, where $\calL$ is some linear differential operator. Assuming $\vx_k$ is one of the nodes in the stencil $P_k$, denote it now as the point $\vx_{\calI^k_j}$, where $\calI^k_j = k$. With this notation, we can now select the appropriate weights from the matrix $\hat{W}^{\phi,k}$ to approximate $\calL$ using the vector ${\bf w}^{k,j}$ defined in Section \ref{sec:stab_leb}. We write the error estimate as
\begin{align}
\lf|\calL f(\vx_{\calI^k_j}) - \sum\limits_{i=1}^n {\bf w}^{k,j}_i f_{\calI^k_i}\rt| \leq \Lambda^k_{\calL} (\vx_{\calI^k_j}) T \max\limits_{|\alpha| = s} \|D^{\alpha} f\|_{W^{1,\infty}(\Omega_k)} \lf( h(\vx_{\calI^k_j})\rt)^s ,
\label{eq:error_estimate}
\end{align}
where $f \in W^{s,\infty}(\Omega_k)$, and $h(\vx) = \max\limits_j \|\vx - \vx_{\calI^k_j}\|_2$. This estimate assumes that we are performing augmented RBF interpolation with polynomial degree $s \geq m$, where $m$ is the order of the polyharmonic spline.  For an estimate on the error of approximating $\calL$ of order $\theta$ with an RBF augmented with a polynomial of degree $s$, we present a slightly modified version of Eq. 19 from~\cite{DavydovMinimal2016}:
\begin{align}
\lf| \calL f(\vx_{\calI^k_j}) - \sum\limits_{i=1}^n {\bf w}^{k,j}_i f_{\calI^k_i} \rt| \leq P(\vx_{\calI^k_j})\max\limits_{|\alpha| = s} \|D^{\alpha} f\|_{W^{1,\infty}(\Omega_k)} \lf( h(\vx_{I^k_j}) \rt)^{s+1-\theta}, 
\end{align}
where $P(\vx)$ is some \emph{growth function} involving $\Lambda^k_{\calL}(\vx)$~\cite{DavydovMinimal2016}. \changeb{This is a \emph{local} estimate depending on the smoothness of $f$ within $\Omega_k$; this emphasizes the local nature of the RBF-FD method and its potential advantages in dealing with functions of limited smoothness}. This estimate has a similar structure to the standard error estimates for polynomial differentiation. While the above estimate requires $\|f\|_{W^{1,\infty}(\Omega_k)}$, this quantity can be converted to $\|f\|_{W^{1,\infty}(\Omega)}$ under some mild conditions if the local norm of the function is not known; for example, see~\cite{SchabackNativeSpaces}. For a proper derivation of these error estimates, see~\cite{DavydovMinimal2016,DavydovSchaback2016}.

In both these error estimates, the local $\calL$-Lebesgue function $\Lambda^k_{\calL}$ determines how the error is distributed as a function of node placement. Clearly, our proposed stabilization approach also has the potential to reduce pointwise approximation errors. It also follows that that the local $\calL$-Lebesgue function could potentially be used to place stencil nodes in a way as to reduce errors in RBF-FD methods in general. We do not explore this approach here. 

\subsection{Estimating speedup with complexity analysis}
\label{sec:orbffd_complex}

In this section, we compare the computational complexity of stable algorithms, augmented RBF-FD and overlapped RBF-FD. We will then use this comparison to estimate the theoretical speedup of our method over augmented RBF-FD as a function of the overlap parameter $\delta$. For the purposes of this analysis, we ignore accuracy considerations. The results section will discuss the tradeoffs between speedup and accuracy.

The standard unaugmented RBF-FD method involves the following operations per stencil: a) one matrix decomposition of cost $O(n^3)$, and b) one back-substitution of cost $O(n^2)$. For $N$ nodes (and therefore $N$ stencils), the total cost is therefore $C_1 = O\lf(N n^2(n + 1)\rt)$. Stable algorithms typically incur a cost of 10--100 times that of the standard unaugmented RBF-FD method. Using the lower end of that estimate, the cost of stable algorithms scales as $C_2 = O\lf(10 N n^2(n+1) \rt)$.

In contrast, the augmented RBF-FD method involves the following operations per stencil: a) one matrix decomposition of cost $O((n+M)^3)$, and b) one back-substitution of cost $O( (n+M)^2)$. Consequently, its total cost is $C_3 = O\lf(N (n+M)^2 (n+M+1) \rt)$. Dropping the $O$ notation, we have the speedup as
\begin{align}
\frac{C_2}{C_3} = \frac{10 n^2(n+1)}{(n+M)^2(n+M+1)}.
\end{align}
The break-even point is obtained by setting $C_2 = C_3$. This gives us the following cubic equation in $M$:
\begin{align}
M^3 + M^2(3n + 1) + Mn(3n + 2) - 9n^2(n+1) = 0.
\end{align}
In general, real-valued solutions $M(n)$ to the cubic equation are well-approximated by $M \approx n + a$, where $a$ is some small integer. The augmented RBF-FD method is therefore as expensive as a stable algorithm only if $M\approx n$. However, in practical scenarios, $M \lesssim \frac{n}{2}$~\cite{FlyerNS,BarnettPHS}. This implies that augmented RBF-FD is always faster than a stable algorithm, though this discussion ignores the relative accuracy per degree of freedom for these two methods.

The overlapped RBF-FD method incurs a greater cost per stencil than the augmented RBF-FD method, but uses $N_{\delta} < N$ stencils. For each stencil, we still have only one matrix decomposition of cost $O((n+M)^3)$. However, we now have more back-substitutions. Assume that $p_1 = p_2 = \hdots = p_{N_{\delta}} = p$, and let $q = \gamma p$ be the number of nodes (per stencil) retained after $\calL$-Lebesgue stabilization ($\gamma \leq 1$). On average, we thus have $q$ back-substitutions per stencil, for a cost $O( q (n+M)^2)$. The total computational cost of our method is therefore given by $C_4 = O\lf(N_{\delta} (n+M)^2 (n+M+q) \rt)$. We drop the $O$-notation and define the \emph{speedup factor} $\eta$ to be
\begin{align}
\eta = \eta(\delta,n,N,N_{\delta}) = C \frac{C_3}{C_4} = \frac{N (n+M+1)}{N_{\delta} (n+M+q)},
\end{align}
where $C$ is a dimension-dependent quantity. We must now estimate $q$, $p$ and $N_{\delta}$. To estimate $p$ in 2D, let $\rho$ be the average stencil width of a stencil with $n$ nodes. Then, we can define the average stencil area to be
\begin{align}
a = \pi \rho^2.
\end{align}
Let the fill distance be $h$. Assuming quasi-uniformity so that the stencil fill distance is equal to the global fill distance, we have $h = \sqrt{\frac{a}{n}}$. Further, given the parameter $\delta$, we can compute the average retention ball area as
\begin{align}
a_{\delta} = \pi r^2 = \pi (1-\delta)^2 \rho^2,
\end{align}
where $r$ is the average retention distance. Assuming that the fill distance in the retention ball equals the fill distance in the entire stencil, we have
\begin{align}
\sqrt{\frac{a_{\delta}}{p}} = \sqrt{\frac{a}{n}},
\end{align}
which gives us the relationship $p = (1 - \delta)^2 n$. In 3D, a similar argument gives us $p = (1 - \delta)^3 n$. Since this formula gives us $p=0$ when $\delta = 1$, we use $p = \max \lf( (1-\delta)^d n, 1\rt)$ to allow us use $\delta = 1$ freely in our descriptions and results.

To estimate $N_{\delta}$, we note that we are given $N$ nodes to distribute among $N_{\delta}$ stencils so that there are only $q = \gamma p$ weights computed per stencil. This implies that $N_{\delta} = \frac{N}{q}$. Using the definition of $N_{\delta}$ in the expression for $\eta$, we have
\begin{align}
\eta &= C\frac{C_3}{C_4} = C \frac{N (n+m+1)}{N_{\delta} (n+M+q)}, \\ 
\implies \eta &= C \frac{q(n+M+1)}{(n+M+q)}.
\end{align}
Setting $p = \max\lf( (1-\delta)^d n, 1\rt)$ in $d$ dimensions, and $q = \gamma p$, we have
\begin{align}
\eta &= C\gamma \frac{\max\lf( (1-\delta)^d n, 1\rt) (n + M + 1)}{n + M + \gamma \max\lf( (1-\delta)^d n, 1\rt)}.
\end{align}
This asymptotic estimate shows that the theoretical speedup is independent of $N$. If we use the overlapping technique in the context of unaugmented RBF-FD, the speedup factor is
\begin{align}
\eta = \hat{C}\gamma \frac{\max\lf( (1-\delta)^d n, 1\rt) (n+1)}{n+\gamma\max\lf( (1-\delta)^d n, 1\rt)},
\end{align}
which is obtained by setting $M=0$ in all our complexity estimates. In general, the value of $\gamma$ is dependent on the exact node distribution, stencil size, and chosen RBF, and is difficult to estimate a priori. For large stencil sizes, stabilization is unnecessary across a wide range of $\delta$ values, implying that $\gamma=1$ in those cases.

\subsection{Estimating the overlap parameter $\delta$}

While $\delta$ can be treated as an input to our method, it can also be estimated from the formula for $p$, the number of nodes one wishes to retain per stencil. Recall that we have $p = \max \lf( (1- \delta)^d n ,1 \rt)$ in $d$ dimensions. Assume $p = t n$, where $t \leq 1$. Dropping the $\max$ notation for convenience, we have a simple polynomial equation: $\lf( 1 - \delta \rt)^d n = t$. Since $\delta$ is always positive, real and never greater than 1, we have $\delta = 1 - \sqrt[\leftroot{-2}\uproot{2}d]{t}$. Thus, if we know the fraction of the nodes we wish to retain per stencil, we can compute $\delta$.  Another possibility is to eliminate $\delta$ from the algorithm entirely, and use the local $\calL$-Lebesgue functions in a greedy algorithm to decide which weights to retain for each stencil. We leave this approach for future work.

\changeb{\section{Time-Stepping}}
\label{sec:mol}
Our focus is on the forced heat equation. The goal is to find a function $c(\vx,t)$ that satisfies the diffusion equation
\begin{align}
\frac{\partial c}{\partial t} = \nu \Delta c + f(\vx,t), \vx \in \Omega,
\label{eq:heat}
\end{align}
where $f(\vx,t)$ is some source term, $\nu$ is the diffusion coefficient, and $\Delta$ is the Laplacian. $c(\vx,t)$ also satisfies boundary conditions on the domain boundary $\partial \Omega$ of the form
\begin{align}
\mathcal{B}c = g(\vx,t), \vx \in \partial \Omega,
\label{eq:heat_bc}
\end{align}
where $\mathcal{B}$ is some (linear) boundary condition operator. For Dirichlet boundary conditions, we have $\mathcal{B} = \mathcal{I}$, the identity operator (signifying point-evaluation); for Neumann boundary conditions, we have $\mathcal{B} = \nu \frac{\partial}{\partial \vn}$, where $\vn$ is the unit outward normal on $\partial \Omega$.

We use a method-of-lines approach to time discretization,\emph{i.e.}, we first approximate the spatial differential operators (including the boundary operators) with overlapped RBF-FD, then solve the resulting set of ordinary differential equations (ODEs) using an implicit time-stepping scheme. For this article, we use the Backward Difference Formula of Order 4 (BDF4)~\cite{Ascher97}. High-order BDF methods require bootstrapping with either lower-order BDF methods, or Runge-Kutta methods. We use the former approach here.

We now describe our method of lines formulation. Without loss of generality, consider discretizing \eqref{eq:heat} and \eqref{eq:heat_bc} using the Backward Euler time integrator. Let $X = [X_i,X_b]^T$ be the set of nodes on the domain, with $X_i$ being the set of interior nodes, and $X_b$ the nodes on the boundary. Further, let $N_i$ be the number of interior points, and $N_b$ the number of boundary points, with the total number of points being $N = N_i + N_b$. We partition the $N \times N$ discrete Laplacian $L$ into the following 4 blocks:
\begin{align}
L = 
\begin{bmatrix}
L_{ii} & L_{ib} \\
L_{bi} & L_{bb}
\end{bmatrix}.
\end{align}
The blocks on the first row correspond to differentiation in the interior, and the blocks on the second row correspond to differentiation on the boundary. However, since we must actually approximate the operator $\mathcal{B}$ on the boundary, we omit blocks $L_{bi}$ and $L_{bb}$. These must be replaced by the $N_b \times N$ block matrix
\begin{align}
B = 
\begin{bmatrix}
B_{bi} & B_{bb}
\end{bmatrix},
\end{align}
where $B$ is obtained by using RBF-FD to approximate $\mathcal{B}$ on the boundary. While the overlapping approach could be used on the boundary, this does not result in a significant speedup relative to overlapping in the interior. We restrict ourselves to the latter, leaving us with a total number of stencils $N_{\delta} = \lf(N_i\rt)_{\delta} + N_b$. 

Letting $C = [C_i,C_b]^T$ be the vector of samples of the solution $c(\vx,t)$ on the node set $X$, $F_i$ be the samples of $f(\vx,t)$ in the interior, and $G_b$ be the samples of $g(\vx,t)$ on the boundary,~\eqref{eq:heat} and \eqref{eq:heat_bc} can be written together as the semi-discrete block equation:
\begin{align}
\begin{bmatrix}
\frac{\partial C_i}{\partial t} \\
G_b
\end{bmatrix} = 
\begin{bmatrix}
\nu L_{ii} & \nu L_{ib} \\
B_{bi} & B_{bb}
\end{bmatrix}
\begin{bmatrix}
C_i \\
C_b
\end{bmatrix} + 
\begin{bmatrix}
F_i \\
{\bf 0}
\end{bmatrix}.
\end{align}
The next step is to discretize in time using backward Euler. Let $t_{m+1} = t_m + \Delta t$, where $\Delta t$ is the time-step and $m$ indexes a time level. Using superscripts for time levels, we have
\begin{align}
\begin{bmatrix}
\frac{C^{m+1}_i - C^m_i}{\Delta t} \\
G_b
\end{bmatrix} = 
\begin{bmatrix}
\nu L_{ii} & \nu L_{ib} \\
B_{bi} & B_{bb}
\end{bmatrix}
\begin{bmatrix}
C^{m+1}_i \\
C^{m+1}_b
\end{bmatrix} + 
\begin{bmatrix}
F^{m+1}_i \\
{\bf 0}
\end{bmatrix},
\end{align}
which assumes that $f(\vx,t)$ and $g(\vx,t)$ are known for all time. If these quantities are known only at specific time levels, they can be extrapolated using previous time levels to the $m+1$ level, resulting in an implicit-explicit (IMEX) scheme~\cite{Ascher97}. Rearranging, we have
\begin{align}
\begin{bmatrix}
I_{ii} - \nu \Delta t L_{ii} & -\nu \Delta t L_{ib} \\
B_{bi} & B_{bb}
\end{bmatrix}
\begin{bmatrix}
C^{m+1}_i \\
C^{m+1}_b
\end{bmatrix}
=
\begin{bmatrix}
C^m_i + \Delta t F^{m+1}_i \\
G^{m+1}_b
\end{bmatrix},
\label{eq:time_step_linsys}
\end{align}
where $I_{ii}$ is the $N_i \times N_i$ identity matrix. This is a general form for arbitrary linear boundary conditions $\mathcal{B}$. For the special case of Dirichlet boundary conditions, we have $B_{bi} = {\bf 0}$ and $B_{bb} = I_{bb}$, the $N_b \times N_b$ identity matrix. Rearranging, this reduces the above system to the $N_i \times N_i$ system:
\begin{align}
\lf(I_{ii} - \nu \Delta t L_{ii} \rt) C^{m+1}_i = C^m_i + \Delta t \lf(F^{m+1}_i + \nu L_{ib} G^{m+1}_b \rt).
\end{align}
For the BDF4 scheme, \eqref{eq:time_step_linsys} becomes:
{\small
\begin{align}
\begin{bmatrix}
I_{ii} - \frac{12}{25} \nu \Delta t L_{ii} & -\frac{12}{25} \nu \Delta t L_{ib} \\
B_{bi} & B_{bb}
\end{bmatrix}
\begin{bmatrix}
C^{m+1}_i \\
C^{m+1}_b
\end{bmatrix}
=
\begin{bmatrix}
f\lf(C^m_i,C^{m-1}_i,C^{m-2}_i,C^{m-3}_i\rt) + \frac{12}{25}\Delta t S^{m+1}_i \\
G^{m+1}_b
\end{bmatrix}
\label{eq:time_step_linsys2},
\end{align}
where
\begin{align*}
f\lf(C^m_i,C^{m-1}_i,C^{m-2}_i,C^{m-3}_i\rt) = \frac{48}{25} C^m_i - \frac{36}{25} C^{m-1}_i + \frac{16}{25} C^{m-2}_i - \frac{3}{25} C^{m-3}_i.
\end{align*}
}
The linear system \eqref{eq:time_step_linsys2} has a solution only if the time-stepping matrix is invertible, or if the right hand side has a zero in the null-space of the time-stepping matrix (applicable to Neumann boundary conditions). In general, the inverse of this matrix is guaranteed to exist if both the upper left block and the Schur complement are non-singular. For details on the invertibility of such \emph{saddle-point matrices}, we refer the reader to~\cite{Benzi2005}. In practice, we have found that the matrix is invertible. We start the BDF4 scheme with a step each of BDF1 (backward Euler), BDF2, and BDF3.

\changeb{\section{Eigenvalue Stability}}
\label{sec:stab}

We now discuss the stability of our method when approximating the Laplacian $\Delta$ and enforcing derivative boundary conditions.

\subsection{Local $\calL$-Lebesgue functions and boundary conditions}
\label{sec:bc_stab}
\begin{figure}[hptb]
\centering
\subfloat[]
{
	\includegraphics[scale=0.5]{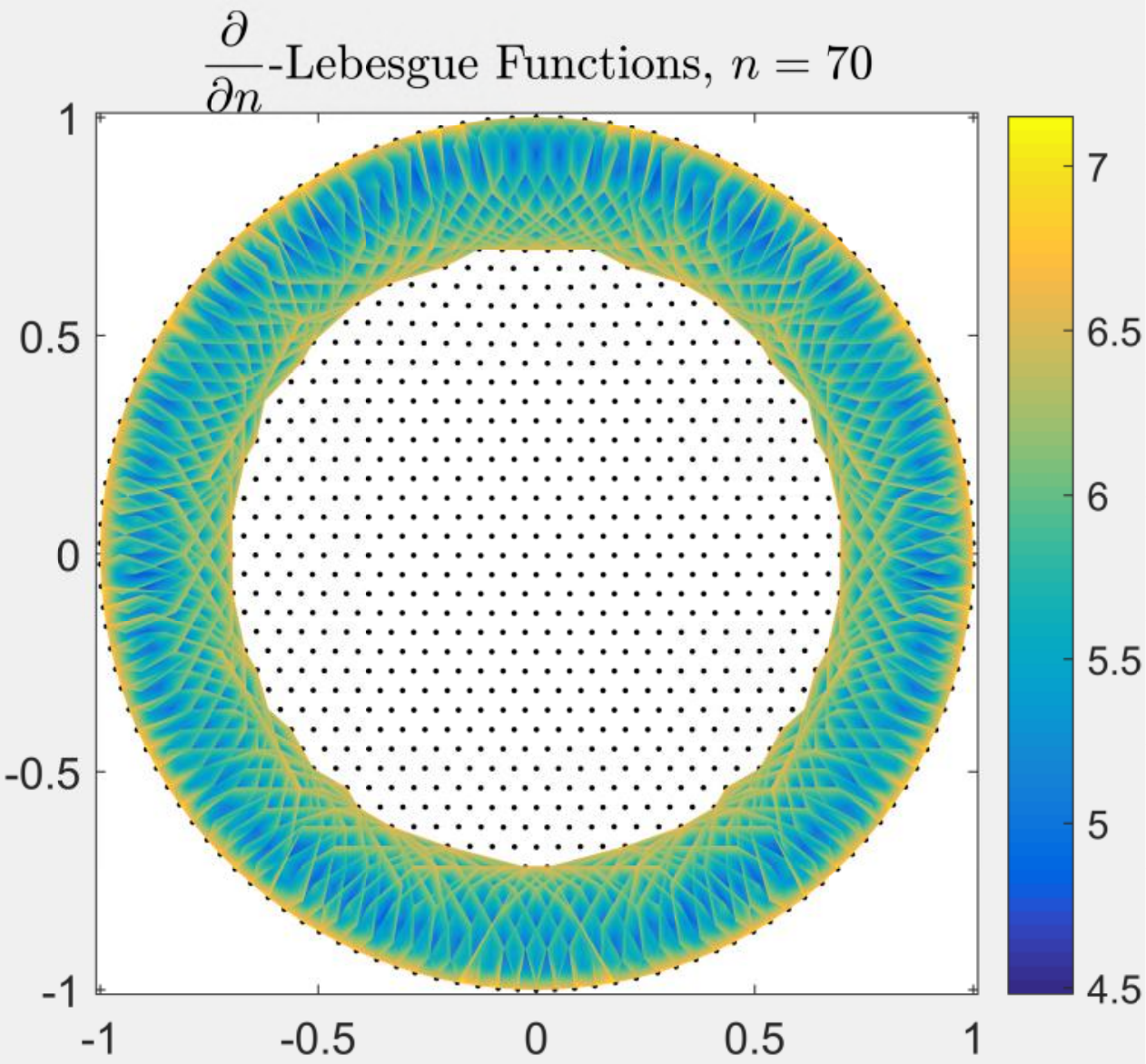} 	
	\label{fig:bc_leba}
}
\subfloat[]
{
	\includegraphics[scale=0.5]{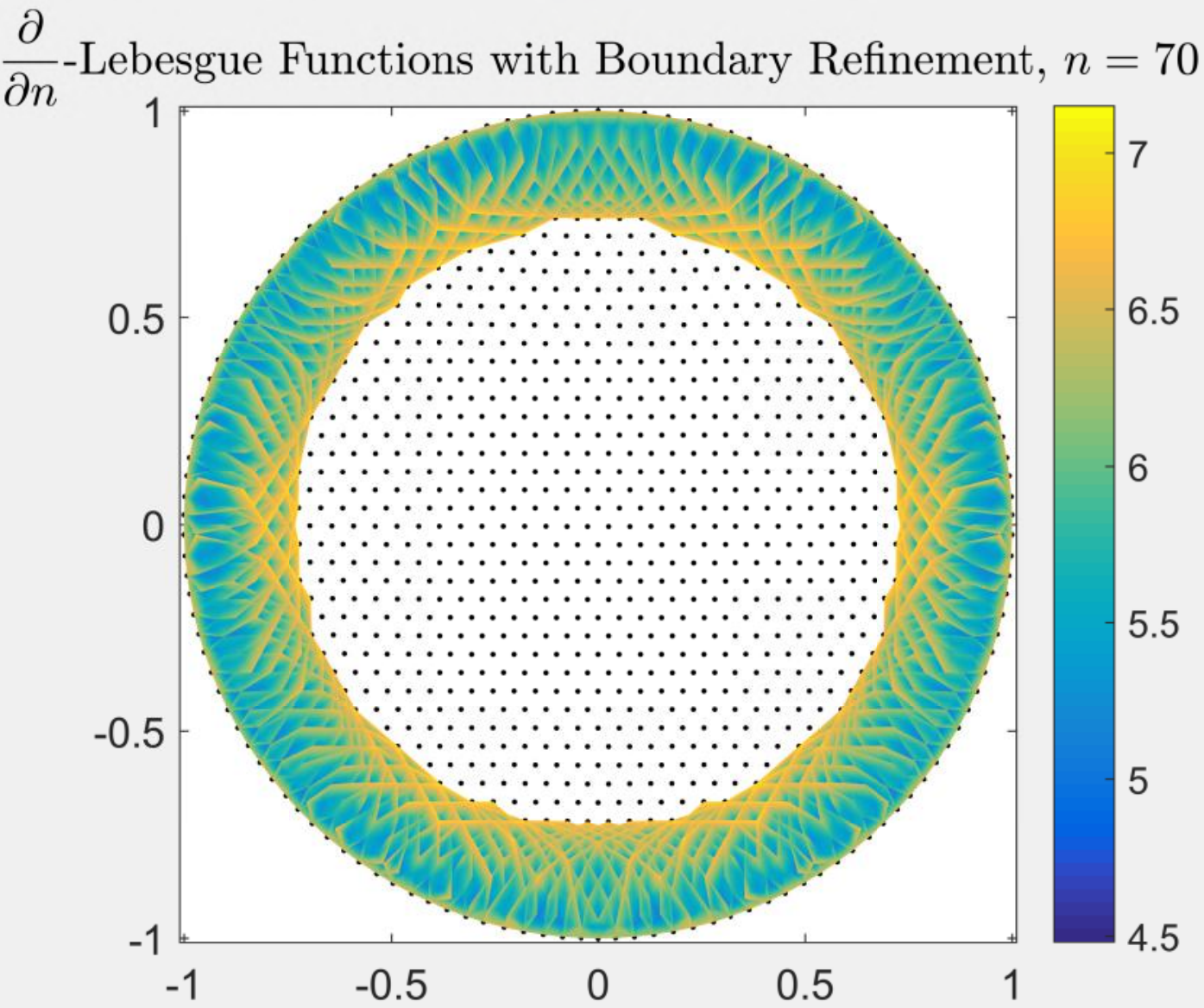}
	
	\label{fig:bc_lebb}
}
\caption{Local Neumann-Lebesgue functions with and without boundary refinement on a logarithmic scale. The figure on the left shows the local $\frac{\partial}{\partial \vn}$-Lebesgue function for boundary nodes on the unit disk. The figure on the right shows the same functions under boundary refinement. Lighter colors indicate higher values.}
\label{fig:bc_leb}	
\end{figure}
The $\calL$-Lebesgue functions give us intuition about stability in the presence of domain boundaries. To better understand this, consider the case of approximating the Neumann operator $\frac{\partial c}{\partial \vn} = \nabla c \cdot \vn$ on the boundary of the unit disk, where $\vn$ is the outward unit  normal. Rather than focusing on approximating a specific function as in Section 3, it should be possible to use the local $\frac{\partial}{\partial \vn}$-Lebesgue functions on each boundary stencil to obtain an intuition for what happens to RBF-FD weights near both stencil and domain boundaries. We set the stencil size to $n=70$, and the total number of nodes to $N=1046$. We visualize the local $\frac{\partial}{\partial \vn}$-Lebesgue functions at every point on each boundary stencil. The results are shown in Figure \ref{fig:bc_leb}.

First, recall from \eqref{eq:error_estimate} that large RBF-FD weights correspond to larger errors through $\Lambda^k_\calL(\vx)$. With this in mind, Figure \ref{fig:bc_leba} confirms that our intuitions from global interpolation in Figure 1 are applicable to RBF-FD: RBF-FD produces large errors on the outermost nodes of each stencil, and progressively smaller errors as one moves towards each the interior of each boundary stencil. This pattern appears in the interior stencils as well (not shown). The presence of dark regions in every stencil around the stencil centers indicates that the overlapping technique is a reasonable generalization of RBF-FD. Unfortunately, on the boundary, the largest weights are exactly the ones we are forced to use when approximating boundary conditions, indicating potential stability issues. 

To alleviate this issue, we add an extra set of points \emph{inside} the domain adjacent to the boundary. This has the effect of \emph{shrinking} the Neumann-Lebesgue function at the boundary nodes; see Figure \ref{fig:bc_lebb}. This approach is known as \emph{boundary refinement}. Figure \ref{fig:bc_lebb} shows the beneficial effect of boundary refinement: notice that the light-colored layer on the boundary has almost vanished. Interestingly, new light-colored regions now appear along the interior edges of the boundary stencils, indicating higher weights there. It is unlikely that this will cause stability issues, since weights in those regions are never computed using boundary stencils. We have noticed spurious eigenvalues in differentiation matrices in the absence of boundary refinement when enforcing Neumann boundary conditions. However, boundary refinement appears to be unnecessary for Dirichlet boundary conditions.

As an alternative to boundary refinement, it is also possible to use ghost points~\cite{FlyerNS}, a set of points outside the domain boundary, to shift the largest weights outside the domain boundary. This is algorithmically more complicated. For simplicity, we do not use this approach for our tests.

\subsection{Local $\Delta$-Lebesgue stabilization}
\begin{figure}[hptb]
\centering
\subfloat[]
{
	\includegraphics[scale=0.5]{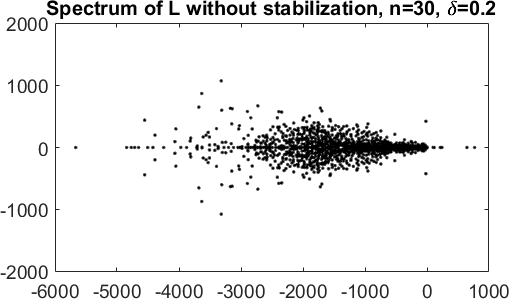} 	
	\label{fig:eig_stab1}
}
\subfloat[]
{
	\includegraphics[scale=0.5]{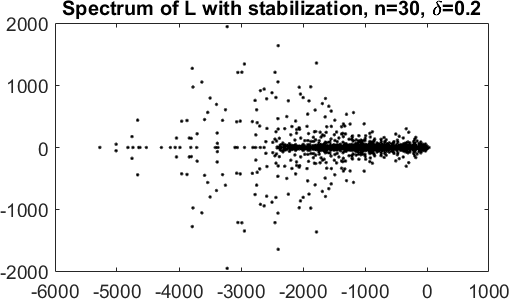}
	
	\label{fig:eig_stab2}
}
\caption{Effect of local Lebesgue stabilization on eigenvalues of the discrete Laplacian. The figure shows shows the spectrum of $L$ without stabilization (left), and with stabilization (right).}
\label{fig:eig_stab}	
\end{figure}
We now present evidence of the efficacy of the local $\calL$-Lebesgue stabilization technique for $\calL \equiv \Delta$ when both $\delta$ and $n$ are small. Again focusing on the disk, we set $n=30$ and use $N=1046$ nodes distributed in the disk. For this value of $n$, it is generally inadvisable to use very small values of $\delta$ lest we retain the large weights seen on stencil boundaries. Indeed, Figure \ref{fig:eig_stab1} shows that $\delta = 0.2$ results in large spurious eigenvalues in this case if no stabilization is used. However, when our automatic stabilization technique is used, the spurious eigenvalue is eliminated (see Figure \ref{fig:eig_stab2}). The trade-off is that the number of stencils has increased, though we still have $N_{\delta} < N$. Stabilization is \emph{not necessary} for larger values of $n$. 
\section{Results}
\label{sec:results}

We apply the overlapped RBF-FD method to the solution of the forced heat equation~\eqref{eq:heat}. The forcing term is selected to maintain a prescribed solution for all time, and the prescribed solution is used to test spatial convergence rates. We solve this test problem using mildly boundary clustered nodes on the closed unit disk in $\mathbb{R}^2$ and the closed unit ball in $\mathbb{R}^3$. These domains were chosen due to the presence of curved boundaries.

\subsection{Forced Diffusion on the disk}
\begin{figure}[htbp]
\centering
\includegraphics[scale=0.6]{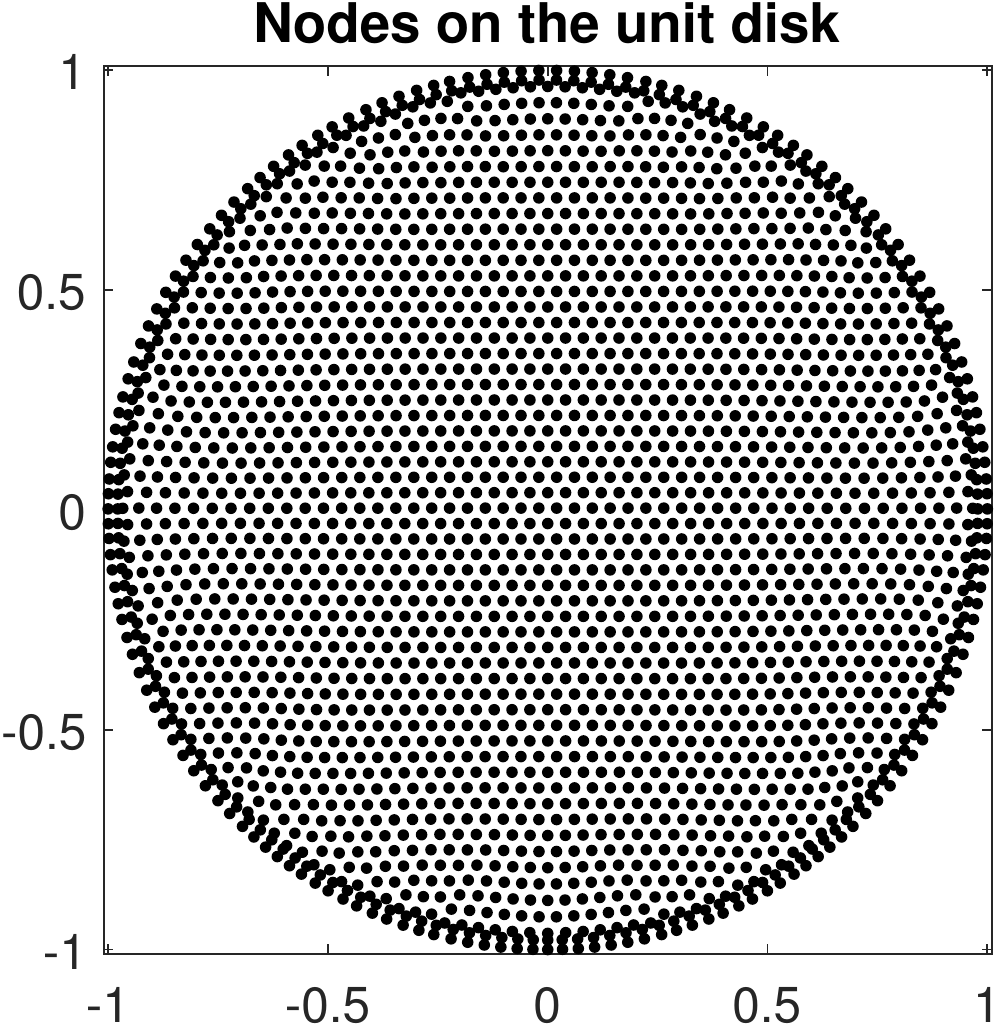} 	
\caption{\changer{Boundary-clustered quasi-hexagonal nodes on the unit disk}. The figure shows the nodes for RBF-FD discretization on the unit disk: interior (solid), near-boundary (cross), and boundary (circles).}
\label{fig:disknodes}	
\end{figure}
We solve \eqref{eq:heat} on the unit disk 
\begin{align}
\Omega(\vx) = \lf\{\vx = (x,y): \|\vx\|_2 \leq 1 \rt\}.
\end{align}
For the following convergence tests, our prescribed exact solution is
\begin{align}
c(\vx,t) = c(x,y,t) = 1 + \sin(\pi x)\cos(\pi y)e^{-\pi t}.
\end{align}
\changer{On the boundary of the domain, we require that $c(\vx,t)$ satisfy an inhomogeneous, time-dependent Neumann boundary condition obtained by evaluating $\nabla c \cdot \vn$ on the boundary $\|\vx\|_2 = 1$, where $\vn$ is the outward normal}. The (interior) forcing term that makes this solution hold is given by
\begin{align}
f(x,y,t) = \pi \lf(2 \pi \nu -1 \rt) \sin(\pi x) \cos(\pi y) e^{-\pi t}.
\end{align}
\changer{To obtain boundary clustered nodes, we generate a node set using the Distmesh program~\cite{DISTMESH}. The resulting nodes are approximately hexagonal in the interior of the domain and have a spacing of approximately $h \propto \frac{1}{\sqrt{N}}$. We then copy the boundary nodes a small distance into the interior along the inward normals $-\vn$ (see Figure \ref{fig:disknodes})}. 

To invert the time-stepping matrix in \eqref{eq:time_step_linsys2}, we use the GMRES iterative method preconditioned with an incomplete LU factorization with zero fill-in, \emph{i.e.,} ILU(0); at each time level, the GMRES method is given the solution from the previous step as a starting guess. With this combination of choices, we have observed that the GMRES scheme typically converges in 10-15 iterations to a relative residual of $O(10^{-12})$ with no restarts required. For the following convergence tests, we set $\nu = 1.0$ and the final time to $T = 0.2$.  The time-step is set to $\Delta t = 10^{-3}$ to ensure that the errors are purely spatial. With this setup, we test accuracy and speedup across a wide range of values of $\delta$, and with $n=30$, $70$, and $101$.

\subsubsection{Accuracy as a function of $\delta$}

\begin{figure}[hptb]
\centering
\subfloat[Convergence rates for $\delta = 1$]
{
	\includegraphics[scale=0.5]{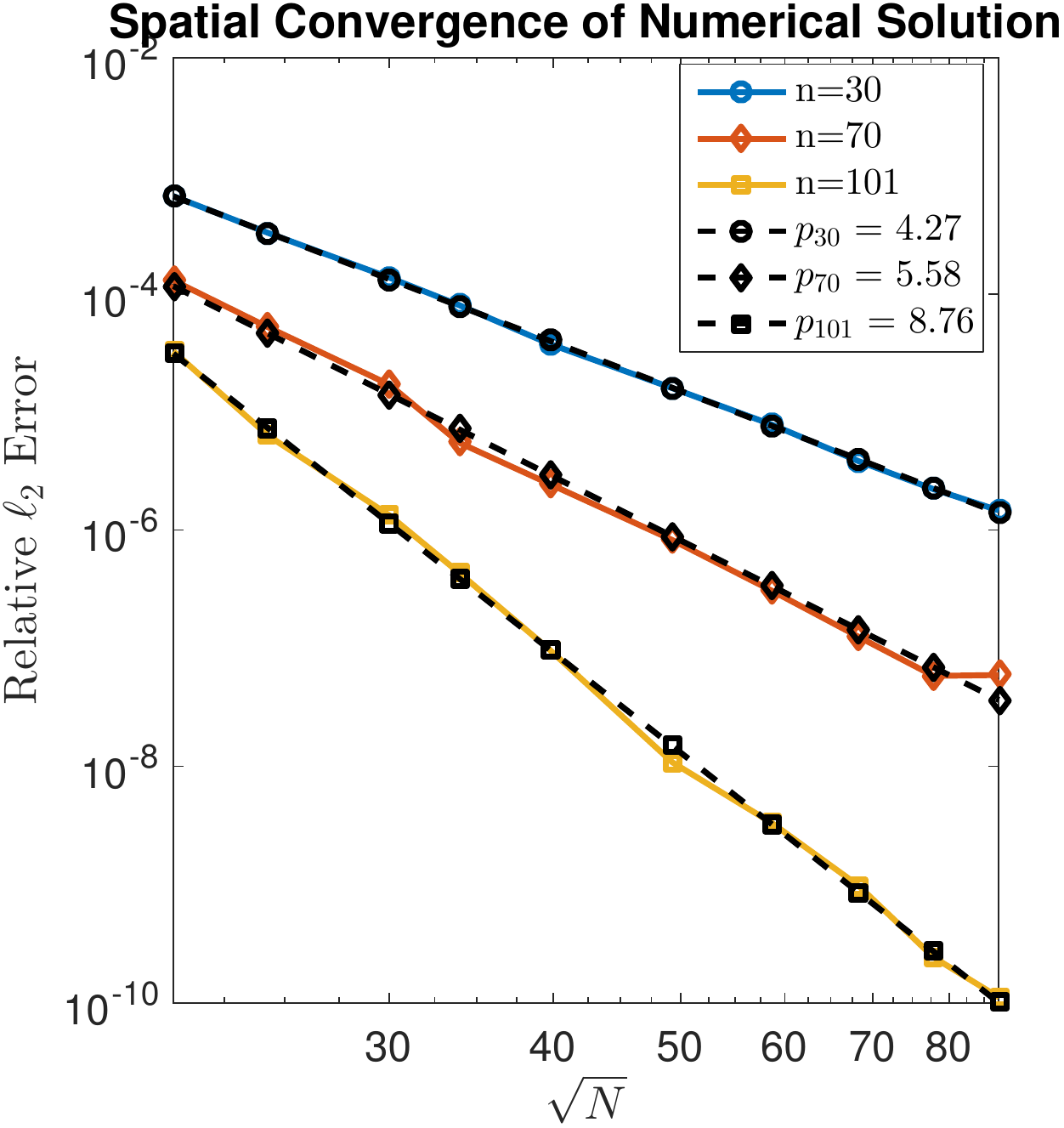} 	
	\label{fig:4a}
}
\subfloat[$n=30$]
{
	\includegraphics[scale=0.5]{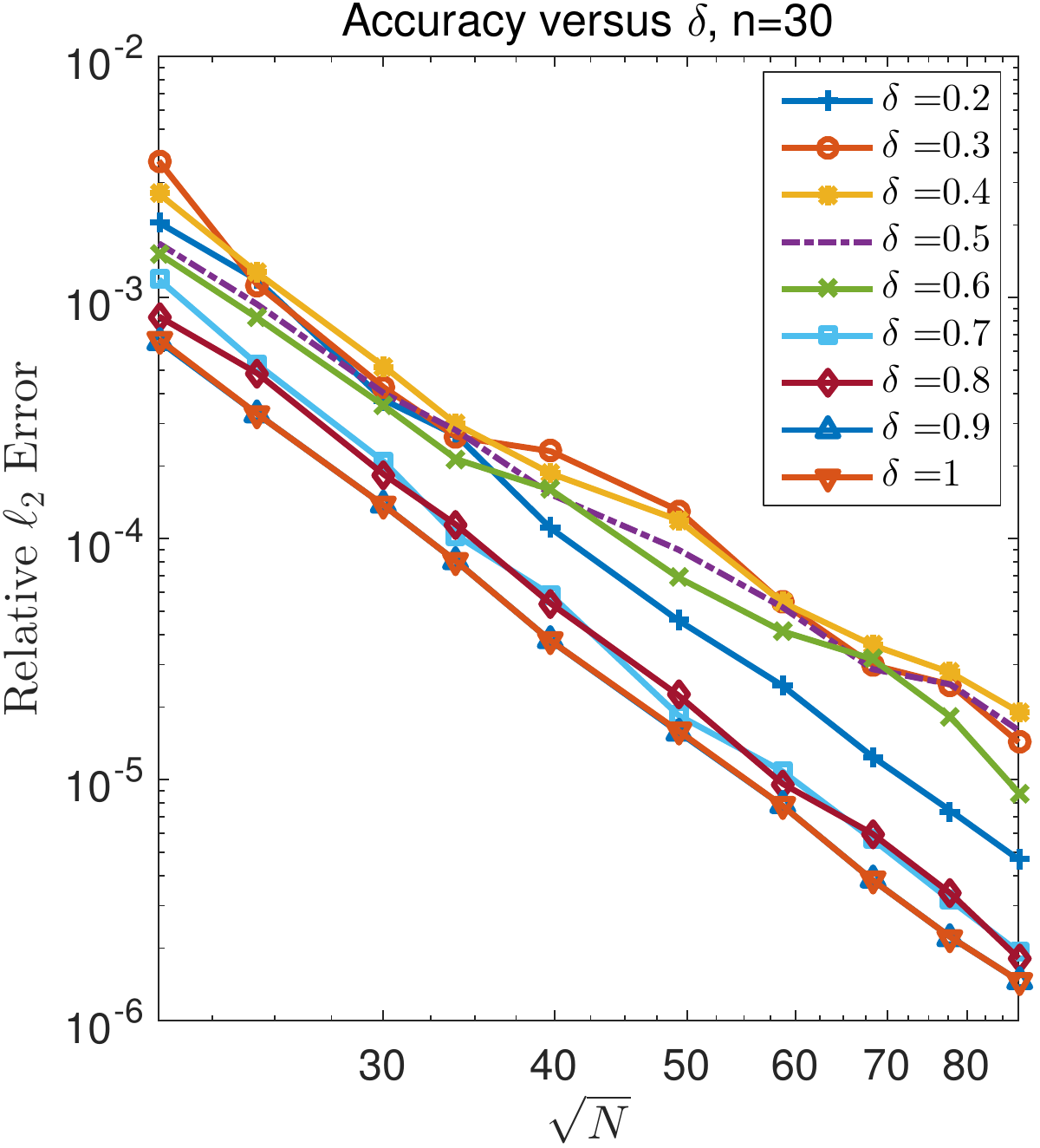}	
	\label{fig:4b}
}

\subfloat[$n=70$]
{
	\includegraphics[scale=0.5]{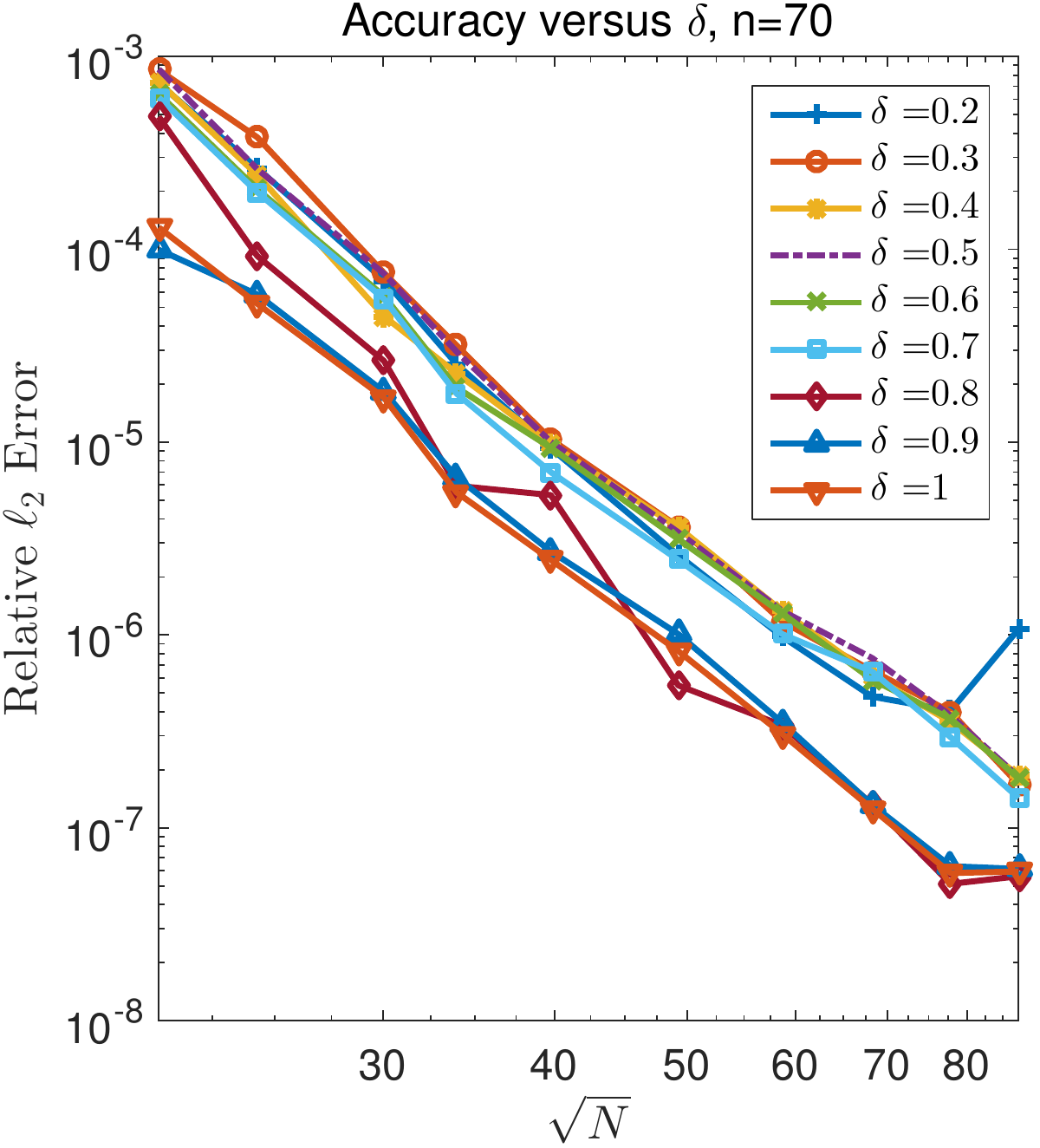}	
	\label{fig:4c}
}
\subfloat[$n=101$]
{
	\includegraphics[scale=0.5]{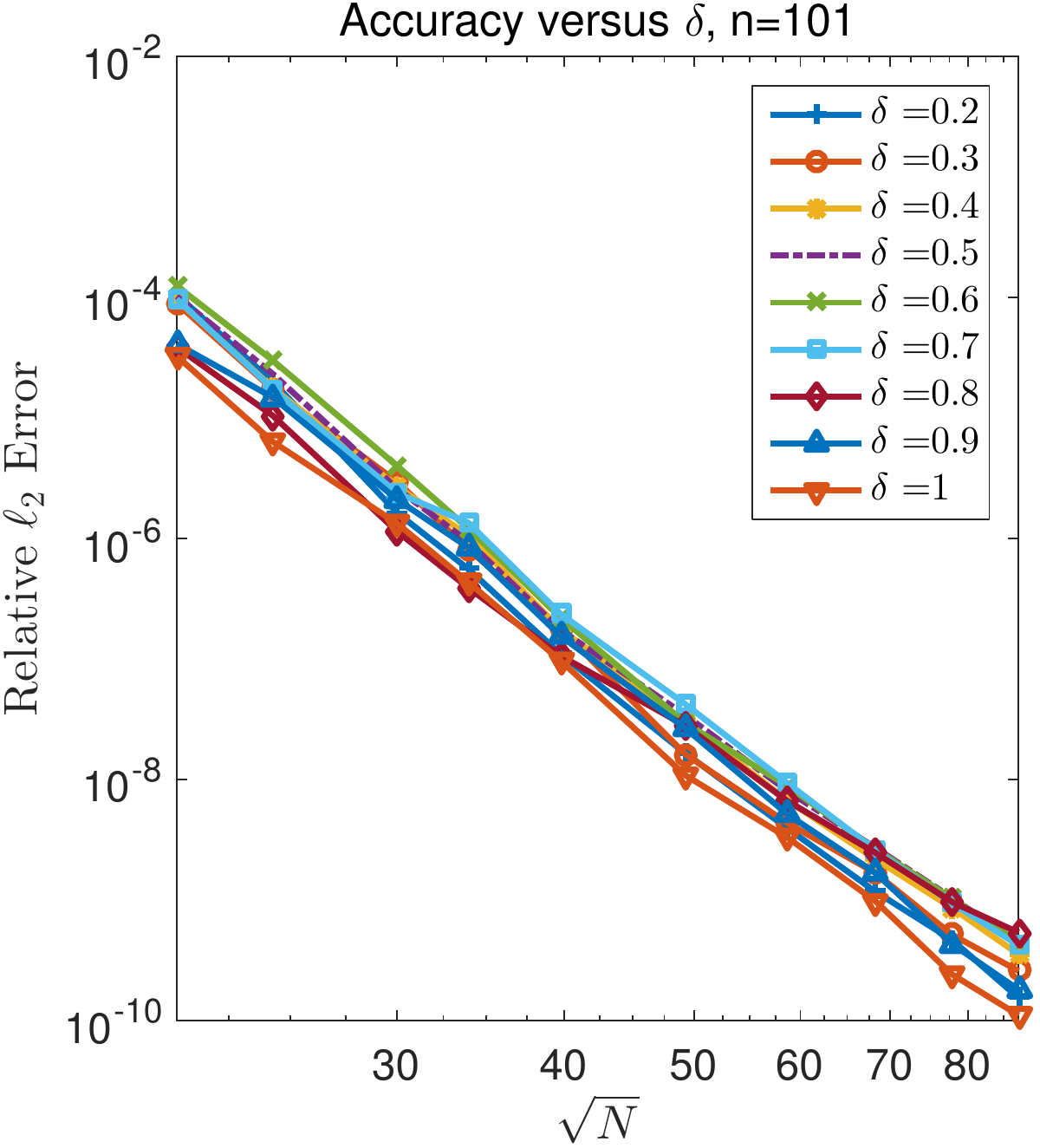}	
	\label{fig:4d}
}
\caption{Accuracy of the overlapped RBF-FD method for the forced diffusion equation on the disk. The figures show $\ell_2$ errors for (a) $\delta = 1$ and $n=30,70,101$; (b) $n=30$, (c) $n=70$, and (d) $n=101$, with $\delta \in [0.2,1]$. The notation $p_n$ indicates the line of best fit to the data for stencil size $n$, showing the approximate rate of convergence for that value of $n$.}
\label{fig:disk_results1}	
\end{figure}

We first test the accuracy of the overlapped RBF-FD method for different values of $\delta$ and $n$. The results of this experiment are shown in Figure \ref{fig:disk_results1}.
Figure \ref{fig:4a} shows different orders of convergence for the RBF-FD method based on the $n$ values used with $\delta = 1$. As the stencil size $n$ is increased, the degree of the appended polynomial increases as well. This increases the order of convergence of the method, giving us third, sixth, and eighth order methods in the $\ell_2$ norm. Similar rates are seen in the $\ell_{\infty}$ norm (not shown). Figures \ref{fig:4b}, \ref{fig:4c}, and \ref{fig:4d} show accuracy as a function of both $\delta$ and the number of nodes $N$ for $n=30$, $70$, and $101$ respectively. When $n=30$, reducing $\delta$ can slightly affect convergence rates, and certainly increase the error by upto an order of magnitude (Figure \ref{fig:4b}). However, only the cases of $\delta = 0.2$ and $\delta = 0.3$ required the use of automatic $\calL$-Lebesgue stabilization. When $n=70$, Figure \ref{fig:4c} shows that convergence \emph{rates} are not affected, and the error increases only slightly. Finally, when $n=101$, Figure \ref{fig:4d} shows that the errors vary very little with $\delta$. 

\subsubsection{Speedup as a function of $\delta$ and $n$}
\label{sec:speed_2d}

\begin{figure}[hptb]
\centering
\subfloat[Theoretical speedup]
{
	\includegraphics[scale=0.5]{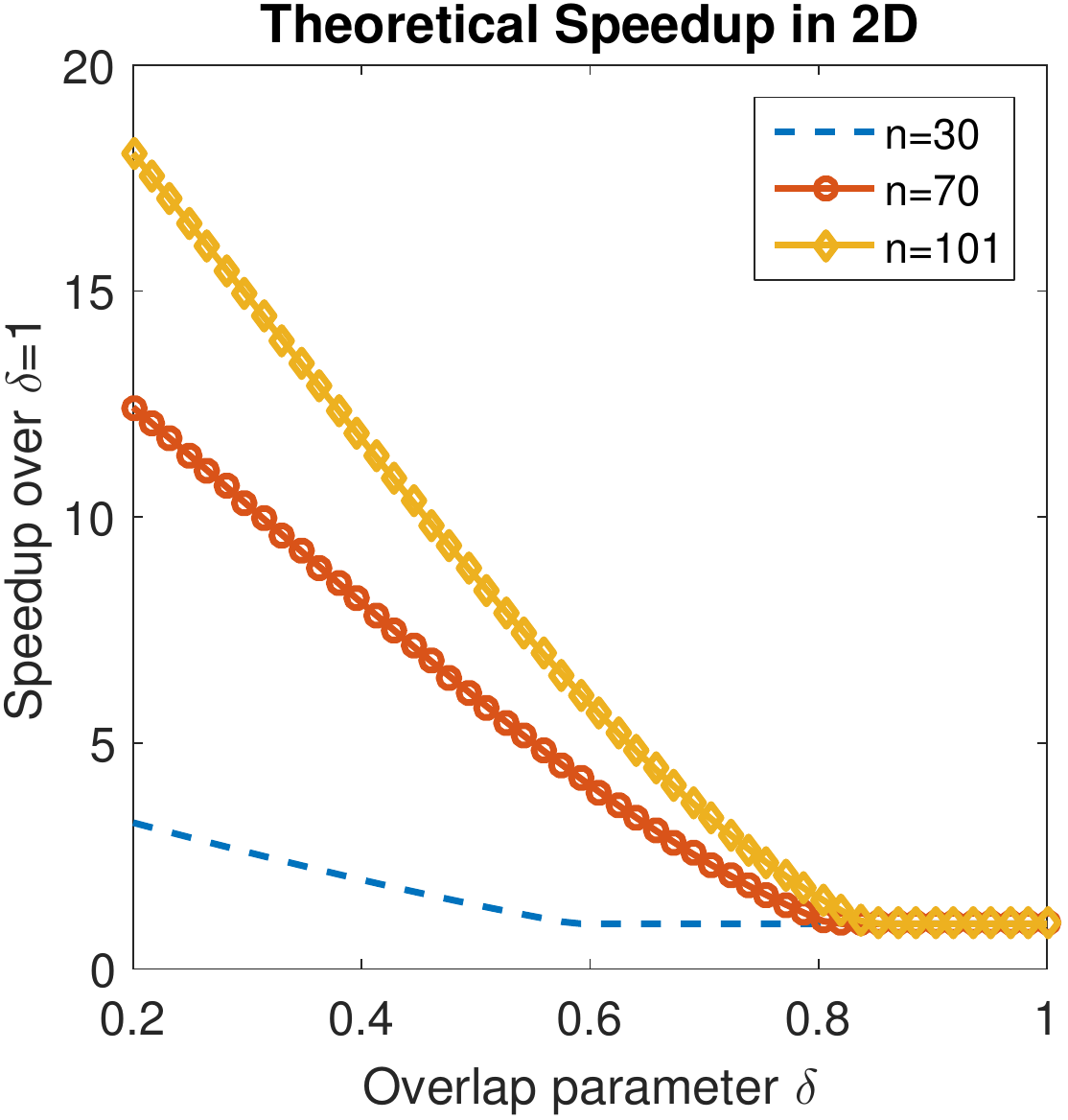} 	
	\label{fig:5a}
}
\subfloat[Observed speedup]
{
	\includegraphics[scale=0.5]{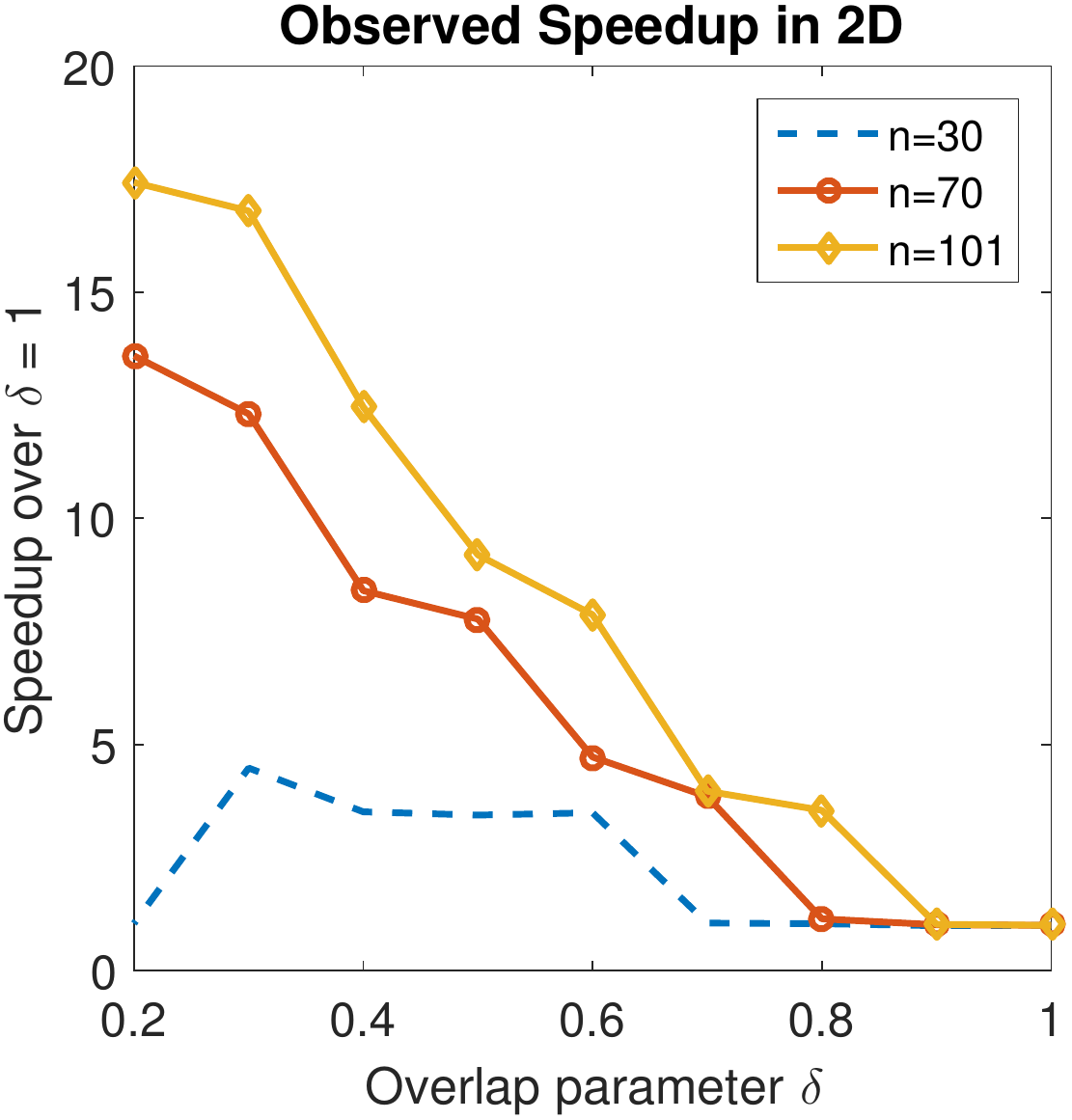}	
	\label{fig:5b}
}

\subfloat[Wall clock time]
{
	\includegraphics[scale=0.5]{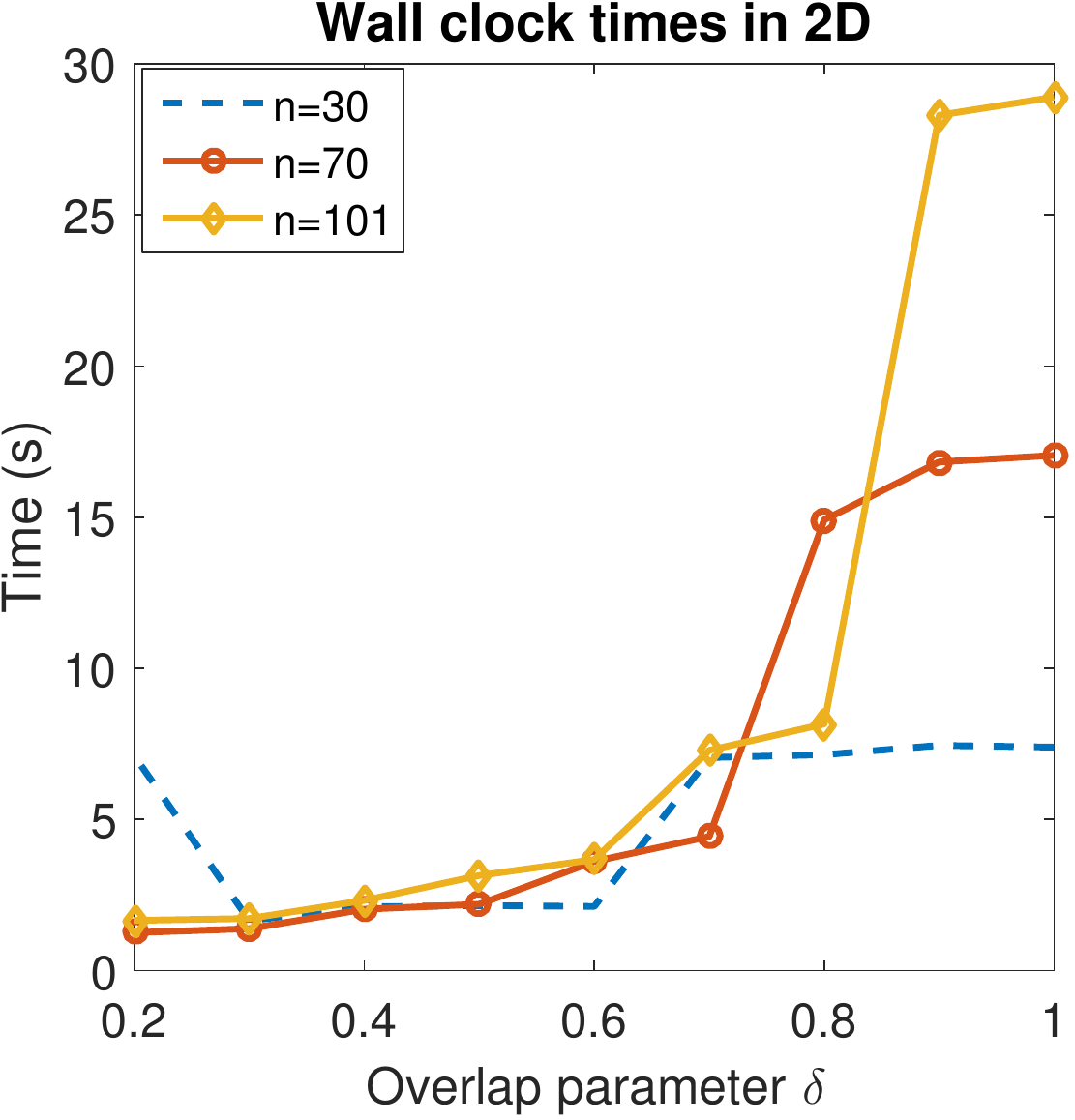}	
	\label{fig:5c}
}
\caption{Speedup of the overlapped RBF-FD method in forming 2D differentiation matrices. The figures show (a) the theoretically estimated speedup, (b) the true speedup, and (c) the actual wall clock time measured in seconds, all as a function of $\delta$ and $n$. Both the theoretical and observed speedups are independent of $N$, the number of nodes.}
\label{fig:disk_results2}	
\end{figure}

We now test the speedup of the overlapped (augmented) RBF-FD method over the augmented RBF-FD method. We compare our theoretical speedup estimate against the speedup obtained in timing experiments. We report speedups and wall clock times so as to not tie our results to language-specific implementations. For $n=30$, we attempt to mimic the effect of stabilization by choosing $\gamma = 0.5$ for our a priori estimate, \emph{i.e.}, $q = 0.5 p$. In all other cases, we set $\gamma = 1$. The constant is set to $C=0.4$ if the theoretical speedup $\eta > 1$. Else, $C=1$. The results are shown in Figure \ref{fig:disk_results2}.
Figure \ref{fig:5a} shows that our a priori estimate of speedup is a good predictor of the observed speedup shown in Figure \ref{fig:5b}. The speedup in \ref{fig:5b} is shown for the highest value of $N = 7615$; however, in practice, we observed that our speedups were independent of $N$ for the values of $N$ tested. MThe maximum speedups are observed for $\delta=0.2$ for $n=70$ and $n=101$. With $n=30$, we notice a drop in speedup, primarily due to rejection of weights from the $\calL$-Lebesgue stabilization algorithm. Combined with Figure \ref{fig:4b}, this indicates that we ought not to use small values of $\delta$ for small values of $n$ in 2D. The highest speedup obtained is approximately 16x that of augmented RBF-FD ($\delta = 1$).

An interesting feature of the overlapped method is that not only are our speedups greater for higher values of $n$, the actual wall clock time is roughly constant as $n$ is increased if $\delta$ is small. This is seen in Figure \ref{fig:5c}. For the case of $\delta = 0.7$, the $n=101$ case is as fast as the $n=30$ case, and far more accurate (\ref{fig:4c},~\ref{fig:4d}), while $n=70$ is actually \emph{less expensive}. Since higher values of $n$ allow for smaller values of $\delta$ without a significant loss in accuracy, it is consequently possible to obtain a higher order method for a cost comparable to that of a low-order method. 
\subsection{Forced heat equation in the unit ball}
\begin{figure}[!ht]
\centering

\includegraphics[scale=0.7]{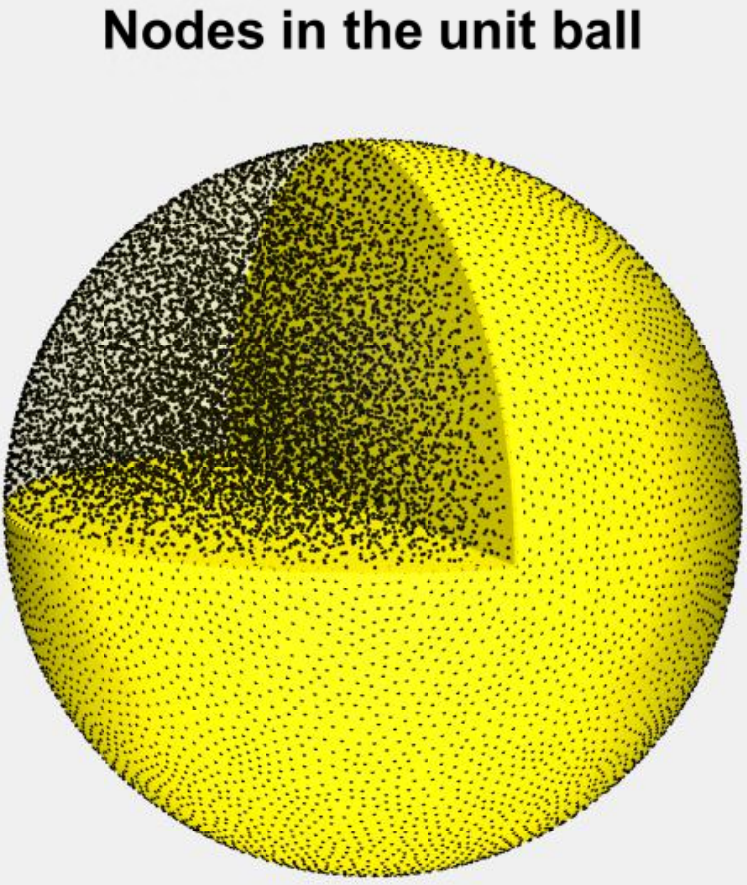} 	
\caption{Nodes in the unit ball. The figure shows the nodes used for the RBF-FD discretization in the unit ball.}
\label{fig:ballnodes}	
\end{figure}
In this test, we solve the forced heat equation in the closed unit ball in $\mathbb{R}^3$:
\begin{align}
\Omega = \lf\{\vx = (x,y,z): \| \vx\|_2 \leq 1 \rt\}.
\end{align}
Once again, we use the method of manufactured solutions. Our prescribed solution is
\begin{align}
c(x,y,z,t) = 1 + \sin(\pi x) \cos(\pi y) \sin(\pi z) e^{-\pi t},
\end{align}
and the corresponding forcing term is
\begin{align}
f(x,y,z,t) = \pi \lf(3\pi \nu - 1 \rt)\sin(\pi x)\cos(\pi y) \sin(\pi z) e^{-\pi t}.
\end{align}
\changer{In this experiment, we use a time-dependent inhomogeneous Dirichlet boundary condition obtained by evaluating $c(\vx,t) = c(x,y,z,t)$ on the boundary of the unit ball, \emph{i.e.,} the sphere $\|\vx\|_2 = 1$}. We obtained node distributions in the unit ball by using the beautiful interactive meshing program Gmsh~\cite{GMSH} to generate a mesh, and then retaining the resulting mesh node distribution. \changer{The resulting irregular nodes have an average spacing of approximately $h \propto \frac{1}{\sqrt[3]{N}}$, and are shown in Figure \ref{fig:ballnodes}. We do not use any boundary clustering in this experiment}.

\begin{figure}[!ht]
\subfloat[Convergence rates for $\delta = 1$]
{
	\includegraphics[scale=0.5]{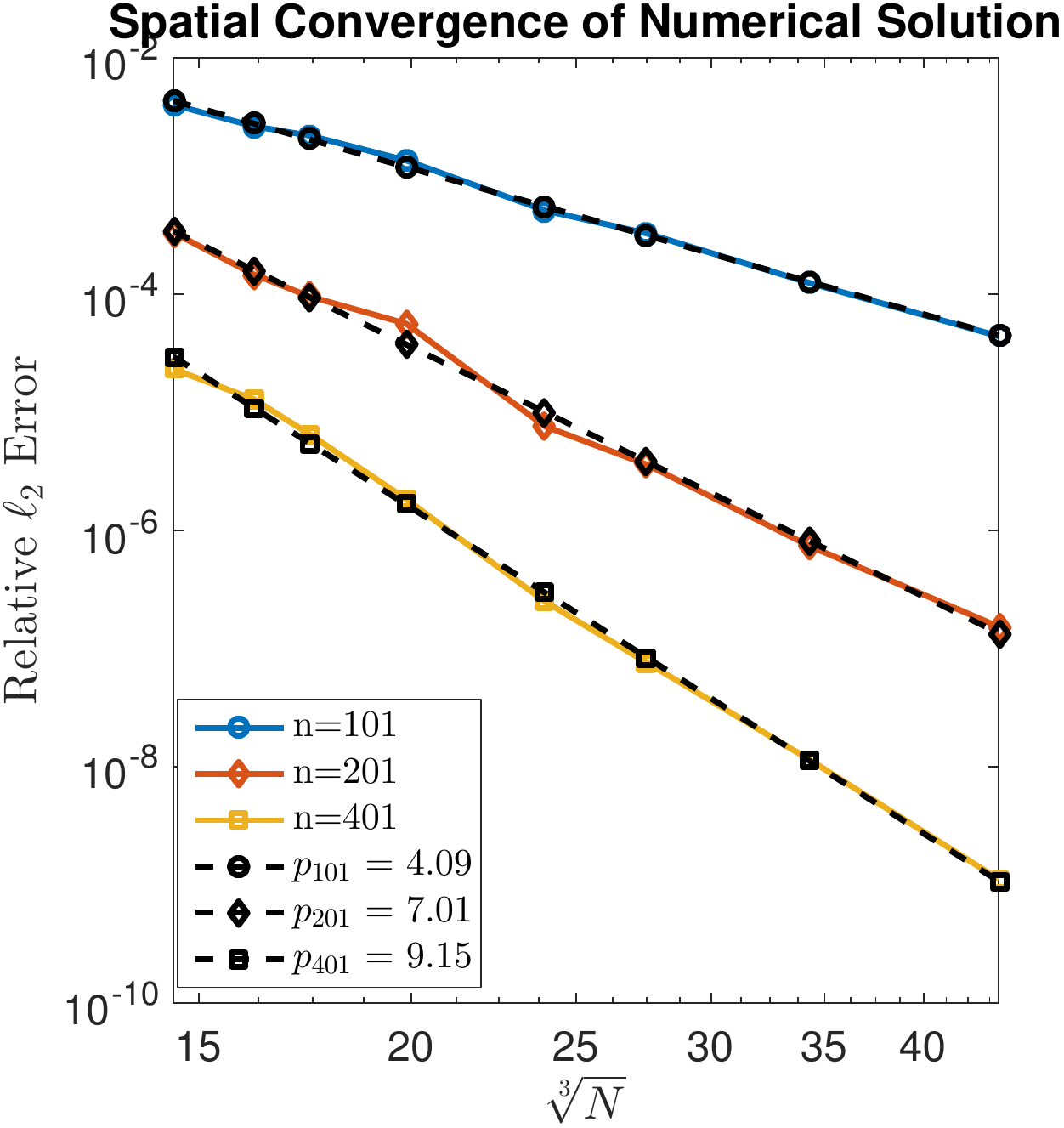} 	
	\label{fig:7a}
}
\subfloat[$n=101$]
{
	\includegraphics[scale=0.5]{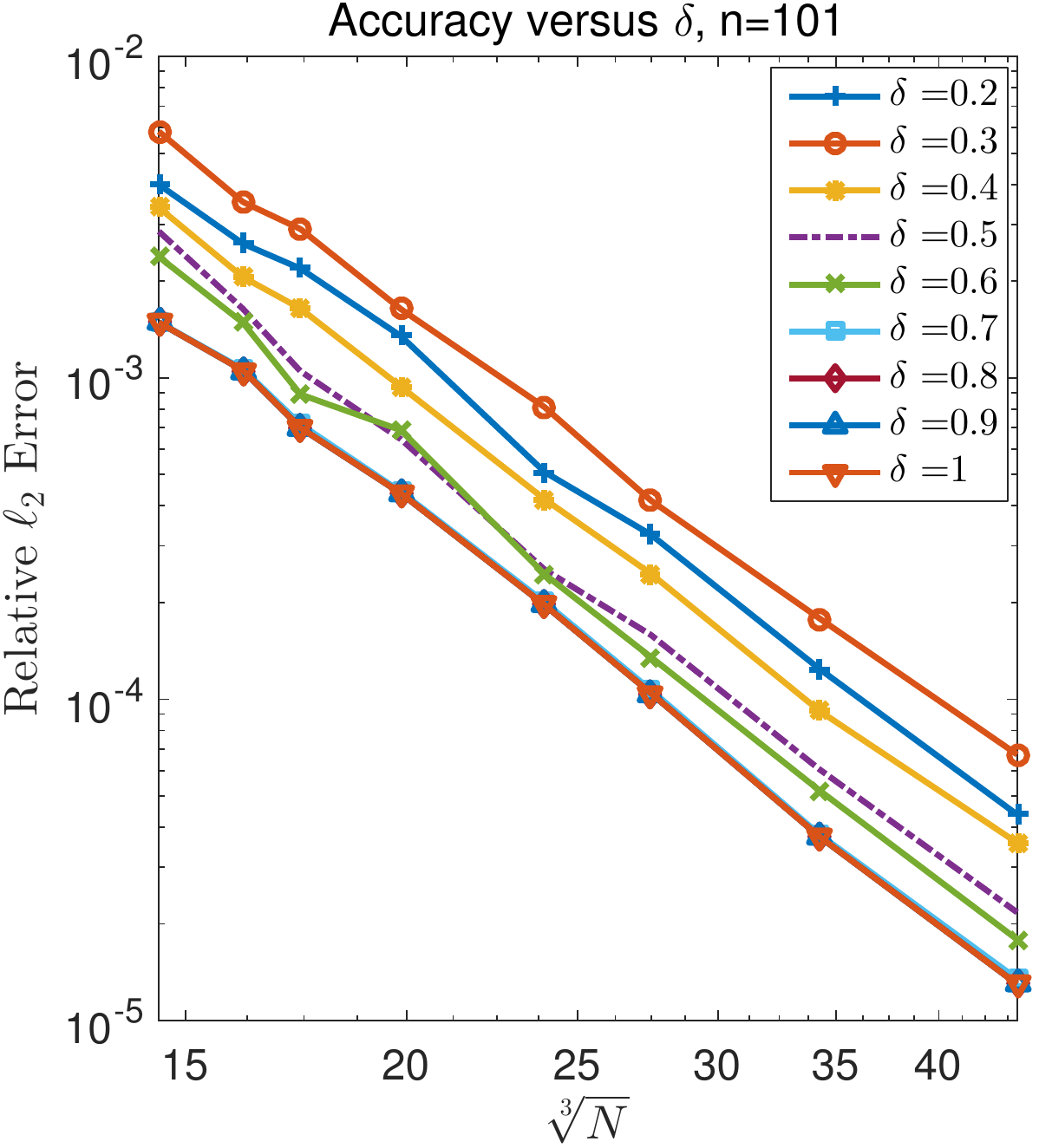}	
	\label{fig:7b}
}

\subfloat[$n=201$]
{
	\includegraphics[scale=0.5]{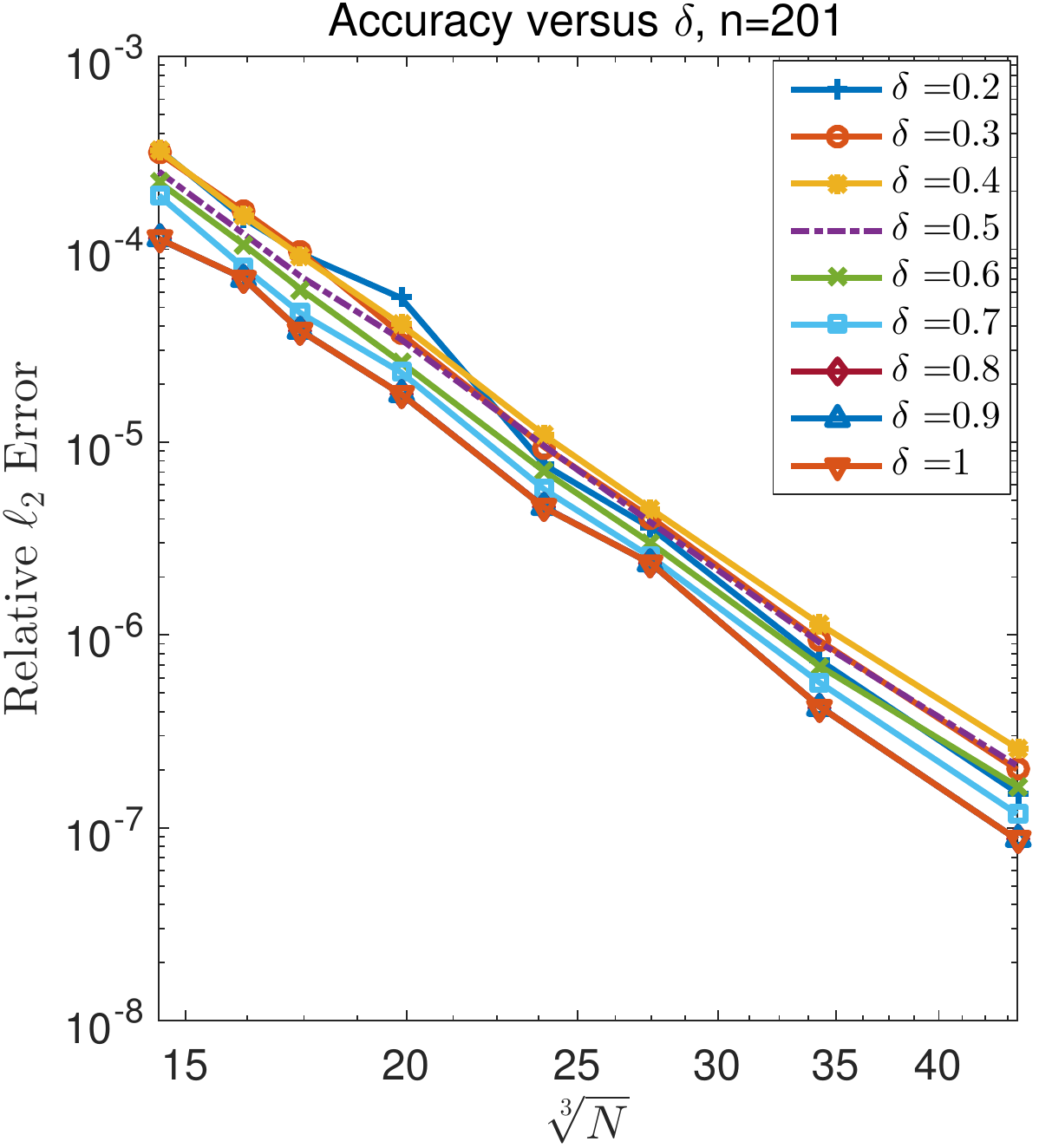}	
	\label{fig:7c}
}
\subfloat[$n=401$]
{
	\includegraphics[scale=0.5]{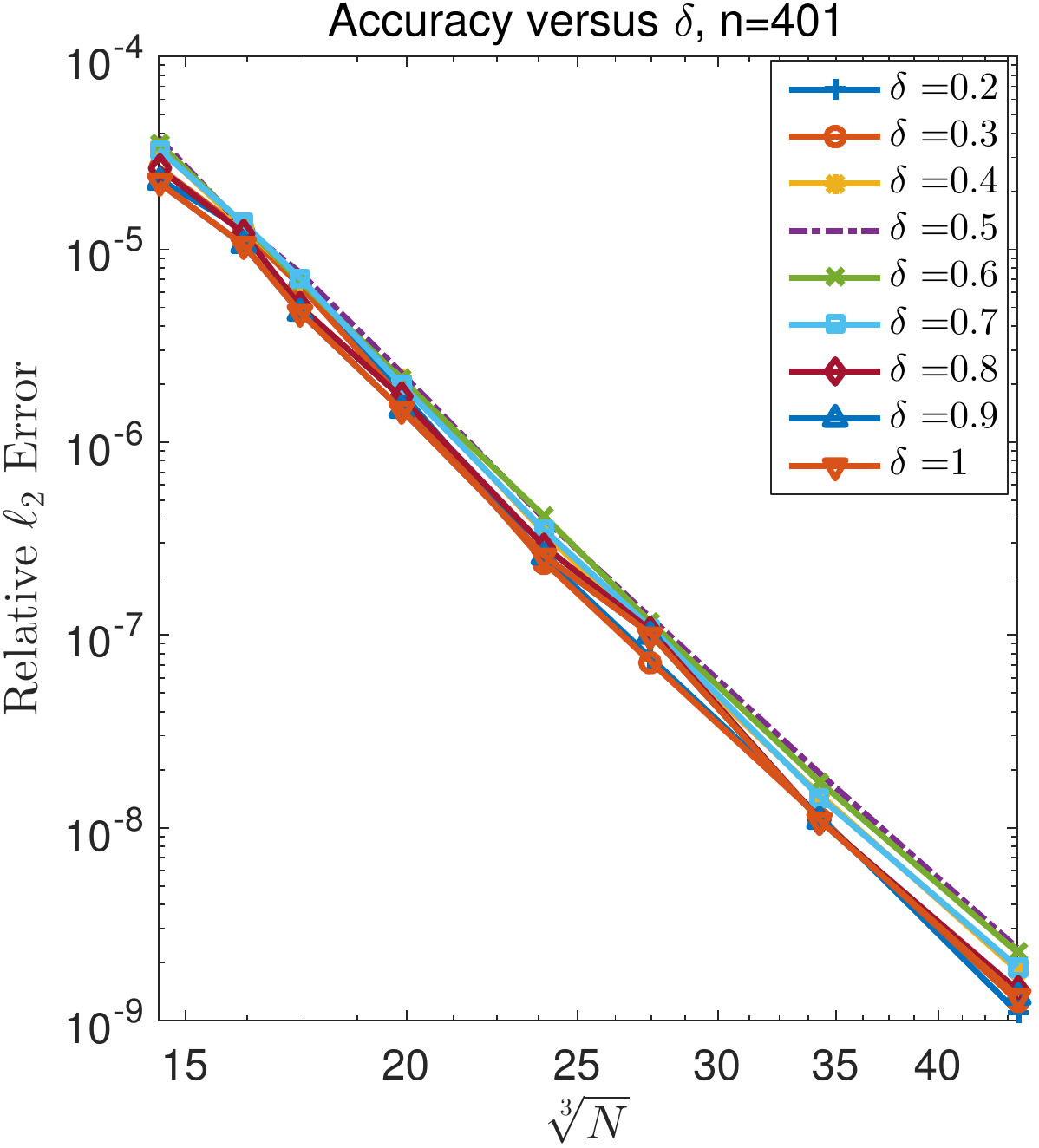}	
	\label{fig:7d}
}
\caption{Accuracy of the overlapped RBF-FD method for the forced diffusion equation in the ball. The figures show $\ell_2$ errors for (a) $\delta = 1$ and $n=101,201,401$; (b) $n=101$, (c) $n=201$, and (d) $n=401$, with $\delta \in [0.2,1]$. The notation $p_n$ indicates the line of best fit to the data for stencil size $n$, showing the approximate rate of convergence for that value of $n$.}
\label{fig:ball_results1}	
\end{figure}

To invert the time-stepping matrix, we now use the BICGSTAB method; this is to avoid the potential memory requirements of the GMRES method in 3D. We precondition the BICGSTAB solver with the ILU(0) factorization of the time-stepping matrix. Further, for each iteration, we feed the solution at the previous time level as a guess. Much like in the case of GMRES, we find that BIGSTAB converges in most cases in 1 or 2 iterations per time-step to a relative residual of $O(10^{-14})$. This significantly lowered iteration count is due to the use of Dirichlet boundary conditions. Once again, we set $\nu = 1$ and the final time to $T=0.2$. The time-step is once again set to $\Delta t = 10^{-3}$. The results are shown in Figure \ref{fig:ball_results1}.
Figure \ref{fig:7a} shows high orders of convergence for large values of $n$. The stencil sizes are much larger than those seen in a 2D problem for the same order. This is due to the fact that we augment the local RBFs with polynomials: if $M \approx \frac{n}{2}$, $n$ must be very large in 3D to support polynomials of even moderate degree. Regardless, we see that it is possible to obtain fourth, sixth, and ninth order methods in the $\ell_2$ norm, with similar results in the $\ell_{\infty}$ norm (not shown). The $n$ values used here are somewhat arbitrary, and other values can be used to obtain similar orders of convergence (as long as the appended polynomial is of the same degree). Figures \ref{fig:7b}, \ref{fig:7c}, and \ref{fig:7d} show accuracy as a function of $\delta$ and $N$ for $n=101$,$201$, and $401$ respectively. Reducing $\delta$ only has a minor effect on the accuracy in 3D. In fact, Figure \ref{fig:7d} shows that the error curves for the different $\delta$ values are clustered very closed together. 

\subsubsection{Speedup as a function of $\delta$ and $n$}
\label{sec:speed_3d}
\begin{figure}[hptb]
\centering
\subfloat[Theoretical speedup]
{
	\includegraphics[scale=0.5]{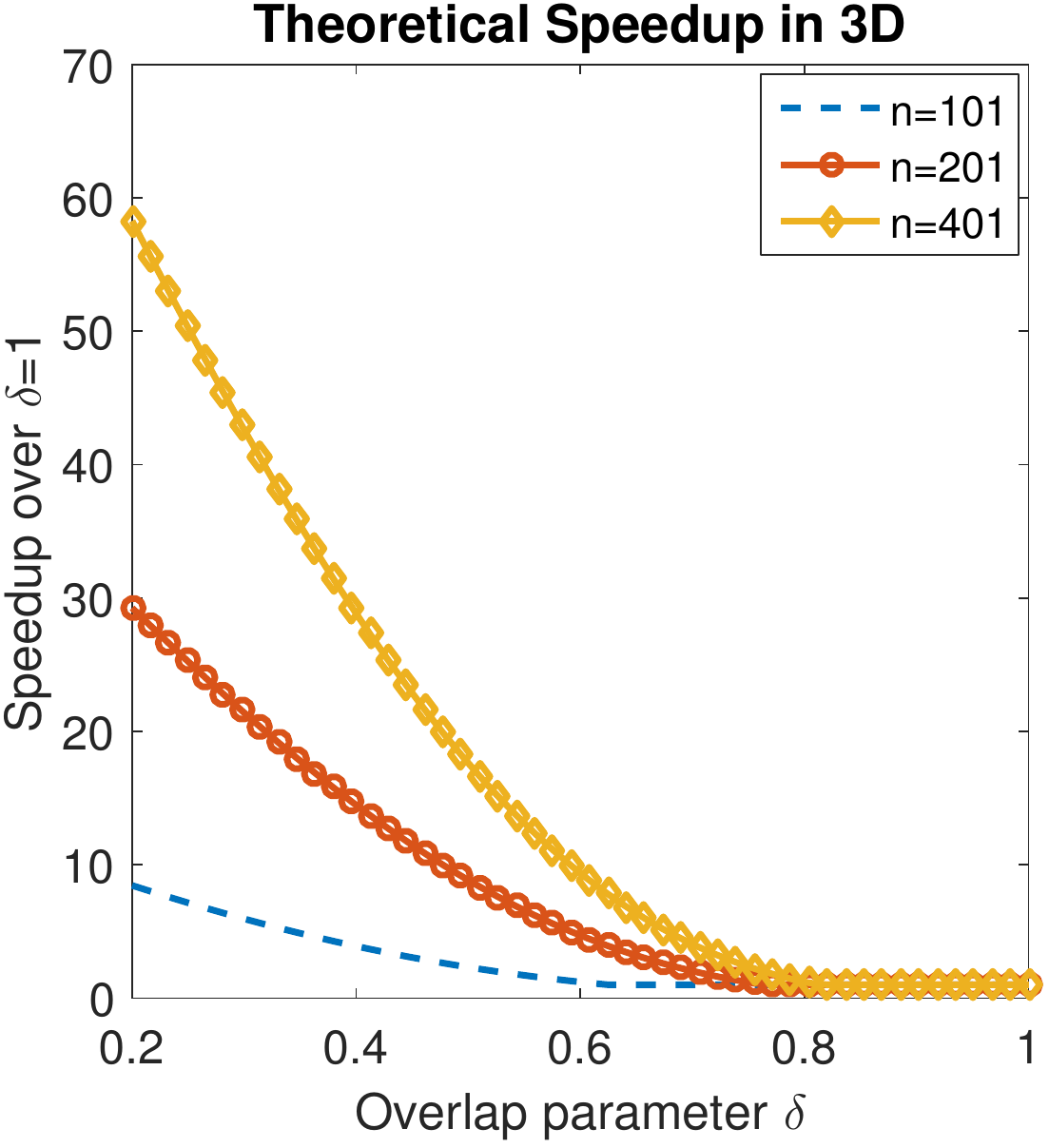} 	
	\label{fig:8a}
}
\subfloat[Observed speedup]
{
	\includegraphics[scale=0.5]{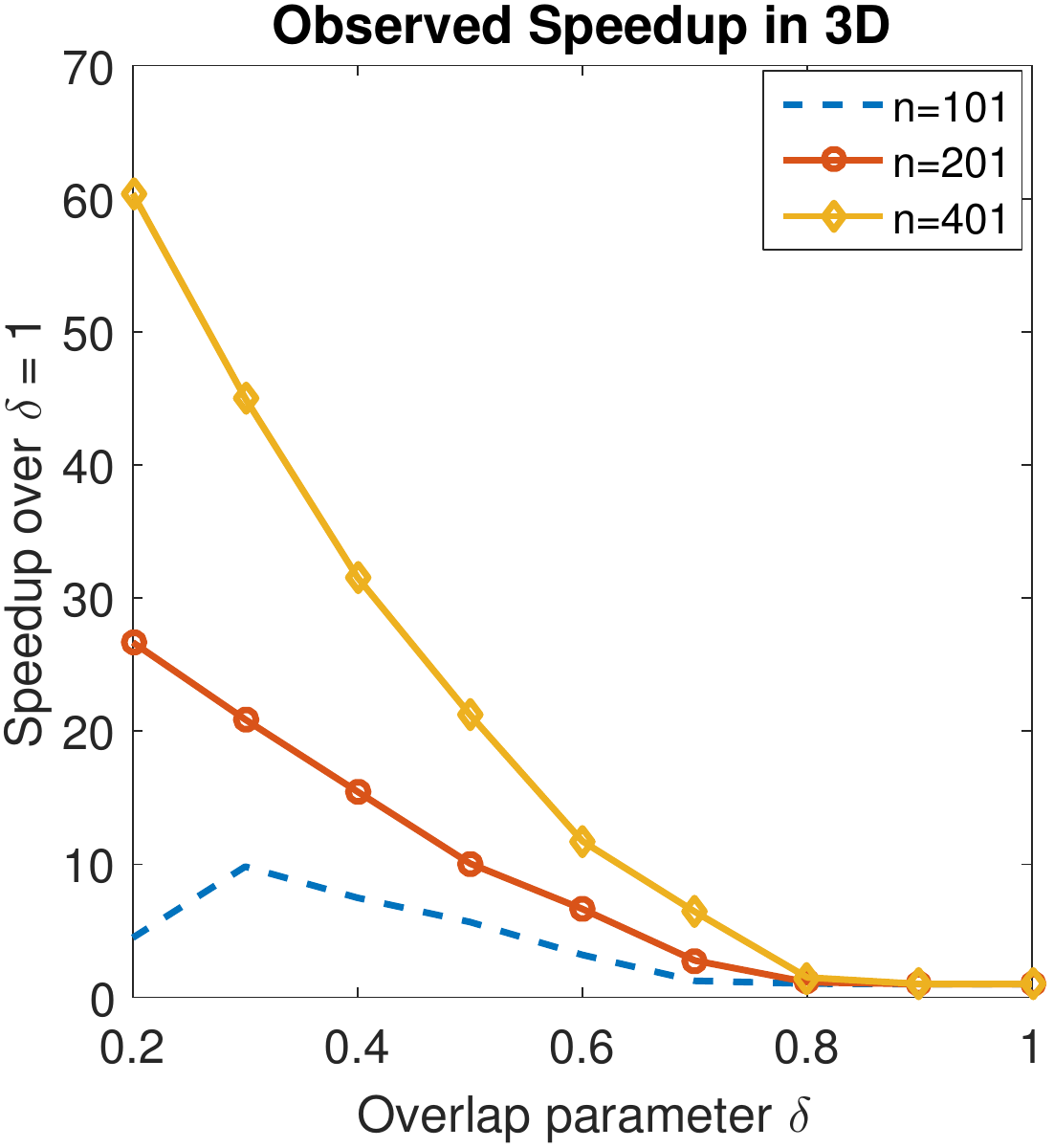}	
	\label{fig:8b}
}

\subfloat[Wall clock time]
{
	\includegraphics[scale=0.5]{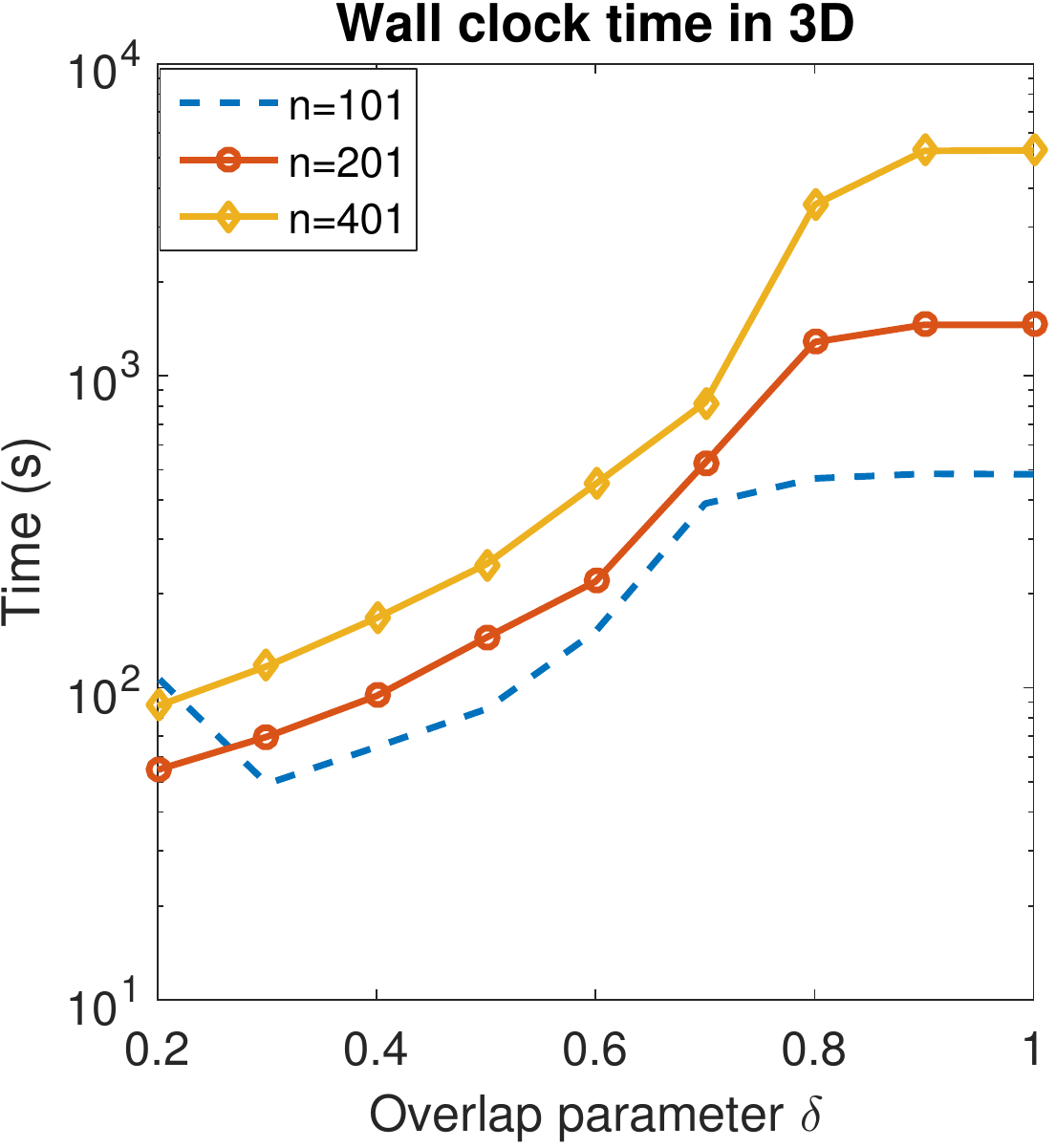}	
	\label{fig:8c}
}
\caption{Speedup of the overlapped RBF-FD method in forming 3D differentiation matrices. The figures show (a) the theoretically estimated speedup, (b) the true speedup, and (c) the actual wall clock time measured in seconds, all as a function of $\delta$ and $n$. Note that both the theoretical and practical speedups are independent of $N$, the number of nodes.}
\label{fig:ball_results2}	
\end{figure}
We now test the speedup of the overlapped (augmented) RBF-FD method over the augmented RBF-FD method. We compare our theoretical speedup estimate against the speedup obtained in timing experiments. We report speedups and wall clock times so as to not tie our results to language-specific implementations. For $n=101$, we attempt to mimic the effect of stabilization by choosing $\gamma = 0.5$ for our theoretical estimate, \emph{i.e.}, $q = 0.5 p$. In all other cases, we set $\gamma = 1$ for our theoretical estimate, as stabilization was not needed. The constant is set to $C=0.38$ if the theoretical speedup $\eta > 1$. Else, $C=1$. The results are shown in Figure \ref{fig:ball_results2}.

Figure \ref{fig:8a} shows that our analytical a priori estimate is once again a reasonable estimate of the observed speedup shown in Figure \ref{fig:8b}. However, the speedups are much higher than in 2D with the greatest speedup being 60x. Again, we observed that our speedups were independent of $N$ for the values of $N$ tested. The maximum speedup is obtained for the highest stencil size of $n=401$, but the smaller stencil sizes also show respectable speedups. The $n=101$ case required stabilization leading to less than optimal speedup.

The actual wall clock time increases very slowly as $n$ is increased, provided that $\delta$ is small (Figure \ref{fig:8c}). For a fixed $N$, the overlapped RBF-FD method reduces the difference between the cost of high-order and low-order methods. Again, since smaller $\delta$ values are feasible for larger $n$ values, it is possible to obtain a high-order method for lower cost than a low-order method; \emph{e.g.}, set $\delta=0.7$ for $n=101$ and $\delta = 0.6$ for $n=201$. In general, the rule of thumb is to decrease $\delta$ as $n$ is increased. The speedup can be predicted in advance using (43).
\section{Summary and Future Work}
\label{sec:summary}

The overlapped RBF-FD method is generalization of the RBF-FD method that helps ameliorate the costs associated with large stencil sizes. It paves the way for very large stencil sizes due to its unique feature of generating high-order methods at a comparable cost to low-order methods. Our method obtained maximum speedups of 16x in 2D and 60x in 3D.

A natural follow-up to our work would be to use the local Lebesgue stabilization technique to compute stencil weights in a greedy fashion on each stencil, thereby eliminating the need for explicitly setting the overlap parameter. It may also be possible to use the local Lebesgue functions to always select the weights in a pattern that enforces conditions on the spectrum of the differentiation matrices. We plan to explore these strategies in a follow-up work. We also plan to explore the relationship between polynomial degree and the polynomial unisolvency of the collocation node set in augmented RBF interpolation. 

The current article only focuses on a serial implementation of the overlapped RBF-FD method. In practical applications, GPU implementations will be necessary. A GPU implementation of the method is currently being developed and compared against a GPU implementation of the standard RBF-FD method.

A natural application of the overlapped RBF-FD method would be the solution of PDEs on time-varying domains and surfaces, where the cost of computing differential operators can no longer be considered a preprocessing step. We are currently exploring the application of our method to a dynamic coupled bulk-surface problem.

\section*{Acknowledgments}
This work was supported by NSF grants DMS-1521748 and DMS-1160432. The author wishes to thank the anonymous reviewers for their detailed suggestions. The author also acknowledges helpful discussions with Professors Akil Narayan (University of Utah), Grady Wright (Boise State University), and Edward Fuselier (High Point University).


\section*{References}
\bibliography{article}

\end{document}